\documentclass[12pt]{amsart}
\usepackage{amssymb}
\usepackage{amsfonts}
\usepackage{amsthm}
\usepackage{amscd}
\numberwithin{equation}{section}
\theoremstyle{plain}
\newtheorem{theorem}{Theorem}[section]
\newtheorem{proposition}[theorem]{Proposition}
\newtheorem{corollary}[theorem]{Corollary}
\newtheorem{lemma}[theorem]{Lemma}
\theoremstyle{definition}
\newtheorem{definition}[theorem]{Definition}

\theoremstyle{remark}
\newtheorem*{remark}{Remark}
\newcommand{\refE}[1]{(\ref{E:#1})}
\newcommand{\refS}[1]{Section~\ref{S:#1}}
\newcommand{\refSS}[1]{Section~\ref{SS:#1}}
\newcommand{\refT}[1]{Theorem~\ref{T:#1}}

\newcommand{\refP}[1]{Proposition~\ref{P:#1}}
\newcommand{\refD}[1]{Definition~\ref{D:#1}}

\newcommand{\refL}[1]{Lemma~\ref{L:#1}}
\renewcommand{\a}{\ensuremath{\alpha}}
\newcommand{\g}{\ensuremath{\gamma}}
\newcommand{\gxy}{\ensuremath{\gamma_{xy}}}
\newcommand{\w}{\ensuremath{\omega}}
\renewcommand{\b}{\ensuremath{\beta}}
\newcommand{\R}{\ensuremath{\mathbb{R}}}
\newcommand{\C}{\ensuremath{\mathbb{C}}}
\newcommand{\N}{\ensuremath{\mathbb{N}}}

\newcommand{\Z}{\ensuremath{\mathbb{Z}}}
\newcommand{\K}{\ensuremath{\mathcal{K}}}

\newcommand{\cint}[1]{\frac 1{2\pi\mathrm{i}}\int_{#1}}
\newcommand{\cintl}[1]{\frac 1{24\pi\mathrm{i}}\int_{#1 }}
\newcommand{\A}{\mathcal{A}}
\newcommand{\hDT}{{h.d.t.}}
\renewcommand{\H}{\mathrm{H}}
\newcommand{\hL}{\mathrm{h.l.}}
\newcommand{\tr}{\mathrm{tr}}
\renewcommand{\i}{\mathrm{i}}
\renewcommand{\l}{\lambda}
\newcommand{\ga}{\mathfrak{g}}
\newcommand{\gb}{\overline{\mathfrak{g}}}
\newcommand{\gh}{\widehat{\mathfrak{g}}}

\newcommand{\KN} {Kri\-che\-ver-Novi\-kov}

\newcommand{\iK}{\ensuremath{1,\ldots, K}}
\newcommand{\Fl}[1][\lambda]{\mathcal{F}^{#1}}
\newcommand{\Fn}[1][\lambda]{\mathcal{F}^{#1}}
\newcommand{\Fln}[1][n]{\mathcal{F}_{#1}^\lambda}

\newcommand{\Ah}{\widehat{\mathcal{A}}}

\renewcommand{\L}{\mathcal{L}}

\newcommand{\Lh}{\widehat{\mathcal{L}}}
\newcommand{\Da}{\mathcal{D}}

\newcommand{\Dalh}{\widehat{\mathcal{D}}_\lambda}
\newcommand{\Do}{\mathcal{D}^1}
\newcommand{\Dh}{\widehat{\mathcal{D}^1}}
\newcommand{\D}{\mathcal{D}^1}
\newcommand{\kndual}[2]{\langle #1,#2\rangle}

\newcommand{\ord}{\operatorname{ord}}
\newcommand{\res}{\operatorname{res}}

\newcommand{\Cal}[1]{\mathcal{#1}}
\newcommand{\fpz}{\frac {d }{dz}}
\newcommand{\ldot}{\,.\,}
\newcommand{\dzl}{\,{dz}^\l}
\newcommand{\pfz}[1]{\frac {d#1}{dz}}
\newcommand{\de}{\delta}
\renewcommand{\d}{\delta}

\renewcommand{\Re}{\mathrm{Re}}
\newcommand{\V}{\mathcal V}
\newcommand{\cinc}[1]{\frac 1{2\pi\mathrm{i}}\int_{#1}}

\newcommand{\sln}{\mathfrak{sl}}

\newcommand{\gl}{\mathfrak{gl}}
\begin{document}
%\baselineskip=20pt

%\layout
%%%%%%%%%%%%%%%%%    private header  %%%%%%%%%%%%%%%%%%%%
\vspace*{-1cm}
\hbox{ }
%\hphantom{\hspace*{\fill} Mannheimer Manuskripte 267}
{\hspace*{\fill} Mannheimer Manuskripte 267}

{\hspace*{\fill} math/0210360}

\vspace*{2cm}

\title[Higher genus affine algebras\quad]
{Higher genus affine algebras 
of Krichever - Novikov type}

\author[M. Schlichenmaier]{Martin Schlichenmaier}
\address[Martin Schlichenmaier]{Department of Mathematics and 
  Computer Science, University of Mannheim, A5, 
         D-68131 Mannheim,
         Germany}
\email{schlichenmaier@math.uni-mannheim.de}
\begin{abstract}
For higher genus multi-point current algebras of Krichever-Novikov type
associated to
a finite-dimen\-sional Lie algebra,  local
Lie algebra two-cocycles are studied.
They yield as central extensions almost-graded
higher genus affine Lie algebras.
In case that the Lie algebra is reductive a complete
classification is given.
For a simple Lie algebra, like in the classical situation,
there is 
up to equivalence and
rescaling only one non-trivial almost-graded central 
extension.
The classification is extended to the algebras of 
meromorphic differential operators of order less or equal  one 
on the currents
algebra.
\end{abstract}
\date{October 17, 2002}
\maketitle 
%%%%%%%%%%%%%%%%%%%%%%%
% section 1
\section{Introduction}\label{S:intro}
%\input intro.tex
%    Introduction     17.10.02
%%%%%%%%%%%%%%%%%%%%%%%%%%%%%%%%%%%%
Current algebras, their central extensions, and their
representations are of great  importance in mathematics,
mathematical physics, and theoretical physics.
These algebras supply examples of 
infinite dimensional Lie algebras, which are still
mathematically accessible.
Typically, they appear in the context of systems with 
infinitely many independent symmetries.
For example, 
for Wess-Zumino-Novikov-Witten theory
current algebras show up as gauge symmetry algebras.
In the process of quantization of these systems 
one is often forced to use  
regularization procedures.
By the ``regularization'' one obtains as quantum 
symmetry algebras  central extensions of the 
original current algebras.

In genus zero conformal field theory, 
based on 
a finite-dimensional
Lie algebra $\ga$,  
(sometimes $\ga$ is called horizontal algebra \cite{FuSchwSym})
the current algebra is 
the algebra of $\ga$-valued Laurent polynomials 
$\gb:=\ga\otimes\C\,[z,z^{-1}]$, with
Lie structure
\begin{equation}
[x\otimes z^n, y\otimes z^m]=[x,y]\otimes z^{n+m},
\qquad x,y\in\ga,\ n,m\in\Z.
\end{equation}
For $\ga$ a simple finite-dimensional Lie algebra
with Cartan-Killing form $(.,.)$ the 
central extension $\gh$ of the current algebra $\gb$ 
is obtained by adding a central element
$t$, i.e. $\gh=\gb\oplus\C\, t$ as vector space,
  and replacing the above structure equation  by
\begin{equation}
[x\otimes z^n, y\otimes z^m]=[x,y]\otimes z^{n+m}-(x,y)\cdot 
n\cdot\delta_{m}^{-n}\cdot t,
\quad [t,\gb]=0,
\quad x,y\in\ga,\ n,m\in\Z.
\end{equation}
This 
central extension $\gh$ is called 
an affine Lie algebra.
In mathematics, e.g. 
they appear in the classification of Kac-Moody algebras
as the untwisted algebras of affine type \cite{Kacorg}, \cite{Kacinf}.
In physics they made their first appearance as 
symmetries of two-dimensional conformal field theory models
for the description of quarks, see
Bardak\c ci and Halpern \cite{BarHal}.
See also, e.g. \cite{FuSchwSym} for further references
in physics.

If $\ga$ is simple, it is well-known
(e.g. see Garland \cite{GarLg}) that $\gh$ is the
unique (up to equivalence and rescaling) nontrivial central
extension  of $\gb$.

The algebras $\gb$ and $\gh$ become graded Lie algebras by
setting the degree to 
\begin{equation}
\deg(x\otimes z^n):=n, \quad \deg t:=0.
\end{equation}
This grading is a crucial property in the construction of representations.

Equivalently, $\gb=\ga\otimes\C\,[z,z^{-1}]$ can 
be described  as the Lie algebra of $\ga$-valued meromorphic functions
on the Riemann sphere (i.e. on the compact Riemann surface of
genus 0) which are holomorphic outside the points $0$ and $\infty$.
Note that the commutative algebra $\C\,[z,z^{-1}]$ 
is the algebra of meromorphic functions on the sphere 
which are holomorphic
outside $0$ and $\infty$.
Starting from this description of the genus zero algebra the
natural extension to higher genus compact Riemann surfaces $M$ 
(or equivalently to smooth projective curves over $\C$)
is to replace $\C\,[z,z^{-1}]$  by the commutative algebra 
$\A$ of those meromorphic functions on $M$ obeying some
regularity conditions.
This was done by Krichever and Novikov
\cite{KNFa, KNFb, KNFc} for the case when the functions have
poles only at two fixed points. 
In more details these algebras were studied by Sheinman
\cite{Shea,Shaff}.
The author of the current article  extended this to 
the case of an arbitrary but finite number of points
where poles are allowed \cite{Schlct,Schlhab}.
In this way we obtain $\gb=\ga\otimes\A$,
where $\A$ is the algebra of meromorphic functions 
which might have poles 
at the chosen set of points.
The Lie algebra is given via
the structure equations
\begin{equation}
[x\otimes f,y\otimes g]=[x,y]\otimes fg\ .
\end{equation}
An essential (and non-trivial) part of these
generalizations is the introduction of a graded structure.
Indeed, in general  it is not possible to introduce a grading, but only an
almost-grading (see \refD{almgrad} below).
But it turns out that for the purpose of representation theory this
is enough.
The almost-grading is introduced by dividing the set $A$ of points where
poles are allowed into two disjoint non-empty subset 
$I$ and $O$.
Different splittings will define non-equivalent almost-gradings.
Clearly, in the two-point case such a splitting and hence
the almost-grading is unique.

As usual central extensions $\gh_\g$ are constructed with the 
help of Lie-algebra 2-cocycles $\g$ for the 
Lie algebra $\gb$ with values in the 
trivial module $\C$.
No other type of cocycles will appear in this article. Hence for
simplicity I will call them just cocycles.

From the point of view of representations we ask for central extensions 
$\gh_\g$  such that we can assign  a degree to  the central element
in such a way that the almost-grading of $\gb$ can be extended to 
$\gh_\g$.
In contrast to the classical situation
(i.e. genus zero and two marked points)
not every cocycle will define a central extension
which allows an extension of the almost-grading.
A necessary and sufficient condition is, that the cocycle is
local in the following sense.
Let $\gb=\ga\otimes\A =\oplus_{n\in\Z}\gb_n$,
  be the decomposition
into the homogeneous subspaces $\gb_n$ with respect to
the almost-grading.
A cocycle $\g$ is called {\it local}, if there exist integers
 $L_1,L_2$, such that
\begin{equation}
\g(\gb_n,\gb_m)\ne 0\implies L_2\le n+m\le L_1.
\end{equation}
By assigning $\deg t:=0$ (or any other integer) the 
almost-grading is extended to $\gh_\g$.
This notion of locality of a cocycle is due to Krichever
and Novikov \cite{KNFa}.

We call  a cohomology class local if it contains a representing
element which is a local cocycle.
Not every representing element 
of a local class will be local.
In the classical case
it turns out, that every cocycle class is local.
This is neither the case if we allow more points nor
if we consider Riemann surfaces of higher genus.
The main goal of this article is to give a complete
classification of local cocycles for the higher genus 
and multi-point current algebras at least in the case when
the finite-dimensional Lie algebra $\ga$ is reductive.
In this way the almost-graded (higher genus and multi-point)
affine algebras are classified as central extensions 
of current algebras up
to equivalence.
Sometimes central extensions defined via local
cocycles are called local central extensions. To
avoid  misinterpretations I prefer to use
{\it almost-graded central extensions} instead.
It should be pointed out that the cocycles
obtained by projective representations are local
in all important cases, e.g. in the
case of semi-infinite wedge representations.

The results obtained in this article are described as follows.
Let $\ga=\ga_0\oplus\ga_1\oplus\cdots\oplus\ga_M$ be a
finite-dimensional reductive Lie algebra with abelian summand $\ga_0$ 
of dimension $n$ and $M$ simple summands $\ga_1,\ldots,\ga_M$.
Let $C_S$ be a cycle separating the points in $I$ from the points in
$O$ (see near Equation \refE{knpair} for  more information
on its definition).
First we consider the semisimple case.
It is shown that  for every local cocycle $\g$ of $\gb$ we can 
find a symmetric invariant bilinear form $\a$ on $\ga$ such that
$\g$ is cohomologous to
\begin{equation}\label{E:ocycl}
\g'(x\otimes f,y\otimes g)=\frac {\a(x,y)}{2\pi\i}\int_{C_S}fdg.
\end{equation}
Vice versa, every such cocycle is local.
The corresponding central extension is equivalent to
the central extension defined by
\begin{equation}\label{E:oext}
[x\otimes f,y\otimes g]:=[x,y]\otimes (fg)+
\frac {\a(x,y)}{2\pi\i}\int_{C_S}fdg.
\end{equation}
This will imply that the space 
$\H^2_{loc}(\gb,\C)$ of local cocycle classes is $M$-dimensional,
see \refT{simple} and
\refT{ss}.
In particular, for the current algebra associated to a simple Lie algebra
 we obtain that there is  up to equivalence and rescaling
 only one non-trivial almost-graded central extension.
In this case $\a$ is necessarily 
the Cartan-Killing form, and is non-degenerate.
Let me again stress the fact, that without the locality assumption
for the cocycle (corresponding to the almost-gradedness of the 
central extension) the statement would be wrong.
The only exception is the classical affine Kac-Moody case.

An important concept shows up, the notion of $\L$-invariance of 
a cocycle.
Here $\L$ denotes the Lie algebra of meromorphic vector fields
on $M$, holomorphic outside the points in $A$.
A cocycle is called {\it $\L$-invariant} if 
\begin{equation}
\g(x\otimes(e\ldot f),y\otimes g)+\g(x\otimes f,y\otimes(e\ldot g))
=0,\quad, \forall  x,y\in\ga,\ f,g\in\A,\ e\in\L.
\end{equation}
Here $e\ldot g$ denotes  the (Lie-)derivative of g
with respect to the vector field $e$.
I will show (see \refT{ss}) that in the semisimple case 
in any local cohomology class there is a unique 
$\L$-invariant representing local cocycle. It will be  given by \refE{ocycl}.
This result might be compared with the notion of
$\mathrm{Diff}\, S^1$-invariant cocycles in the $C^\infty$-context
 studied by Pressley and Segal \cite{PS}.

In the semisimple case, 
restricting the attention to $\L$-invariant cocycles does not reduce the
space of cocycle classes.
This is different in  the general reductive case. Already in the case
of pure abelian $\ga$, e.g. for $\gb=\A$, we obtain many
linearly independent local cocycles.
Because $\gb$ is abelian they cannot be cohomologous.
But if we assume  for the cocycle under consideration that
its restriction to $\gb_0\times\gb_0$ is $\L$-invariant, we
obtain the same form of the classification result as 
in the semisimple case.
Again a symmetric invariant bilinear form $\a$ 
classifies those local cohomology classes,
and the central extensions are equivalent to central extensions
given by 
\refE{oext}.
The condition of invariance for $\a$ is restricted to 
$\ga_0$ void. Hence,
the cohomology space  $\H^2_{loc,\L}(\gb,\C)$ 
has  dimension
\begin{equation}
\frac {n(n+1)}{2}+M,
\qquad   n=\dim\ga_0,\quad M \ \text{is the number of simple summands.}
\end{equation}
For non-reductive  $\ga$ certain results are obtained.
But there is no
complete classification.
 
In the context of fermionic or semi-infinite wedge representations
of $\gb$ , resp. $\gh$, there is always a representation of 
a central extension of the  differential operator algebra $\D_\ga$ 
of operators of degree $\le 1$ around.
The algebra $\D_\ga$ is a semi-direct sum of
$\gb$ and $\L$ induced by the action of
$\L$ on $\gb$,
i.e. $\Do_\ga=\gb\oplus\L$
as vector space with  
given Lie structure on the summands, and $[e,x\otimes f]=x\otimes (e\ldot f)$
for $e\in\L,x\in\ga,f\in \A$.
For $\ga=\C$ this algebra is the standard  algebra of meromorphic
differential operators of degree $\le 1$ with poles at
most at the points in $A$.
For this algebra the local cocycle classes have been determined  in 
\cite{SchlCo}.
The space of local cocycle classes is 3-dimensional. It is 
generated by a  standard cocycle  $\g_S^{(f)}$ for the function algebra
$\A$, a standard cocycle  $\g_S^{(v)}$ for the vector field algebra $\L$, and
a mixing cocycle $\g_S^{(m)}$.
The results are recalled in \refS{kn} from \cite{SchlCo}.
In the classical case the corresponding result was
shown by  Arbarello, De~Concini, Kac, and Procesi \cite{ADKP}.

Using these results and the classification result 
for the cocycles of $\gb$ for 
$\ga$ reductive,
with some additional work, a complete classification 
of local cocycles of $\D_\ga$ is given.
It turns out, that every local cocycle $\g$ is cohomologous to
a linear combination of a local (and $\L$-invariant)
cocycle for $\gb$, of the   vector field  cocycle $\g_S^{(v)}$, and a cocycle 
 $\g_\phi$ of mixing type, defined by
\begin{equation}
\g_\phi(e,x\otimes g)=\phi(x)\cdot\g_S^{(m)}(e,g), \quad e\in\L,g\in A,x\in\ga.
\end{equation}
Here $\phi$ is a linear form on $\ga$ vanishing on 
the derived subalgebra
$\ga'=[\ga,\ga]$, i.e. on  the semisimple complement of $\ga_0$.
It is understood that all these basic cocycles are extended by zero
on the corresponding complementary subspaces.
It should be noted that the condition of $\L$-invariance of $\g$ 
restricted to $\gb_0$ is automatically fulfilled, as
$\g$ is a cocycle for  $\Do_\ga$, see \refP{linvg}.
Note also, that if $\ga$ is a semisimple Lie algebra the 
linear form $\phi$ has to vanish identically, hence there will be no 
cohomologically non-trivial cocycles of 
mixing type.
The absence of the mixing cocycle in the semisimple case was also
shown by Sheinman in \cite{Shsc}.

Such central extensions $\widehat{\D_\ga}$ of $\D_\ga$ show also up in the
context of the Sugawara representations associated to 
representations of higher genus affine algebras
\cite{SchlShSu}, \cite{Schli:SOHGRS}.

Algebras of the type $\gb=\ga\otimes \A$ for an arbitrary commutative algebra 
$A$ and  $\ga$ a finite-dimensional simple Lie algebra, 
and their central extensions were considered earlier from the
purely algebraic point of view 
by  Kassel \cite{Kaskd}, and Kassel and Loday \cite{KasLod}.
In particular, no grading was considered.
Hence the question, which cocycles are local and define  
(almost)-graded central extensions could also not be considered.
In \cite{Kaskd} for $\ga$ simple with the help of the 
K\"ahler differentials of $\A$ the universal central 
extension of $\gb$ was constructed.
Indeed, I use this universal central extension 
in the proof of \refT{simple}.
It should be remarked that Bremner \cite{Bremce, Brem4a} gave an 
explicit description of the universal central extensions
in the genus zero case with four marked
points, and in the genus one case with two marked points.
For completeness let me note that  by
Santharoubane \cite{SanSc}, Haddi \cite{HadHo},
Zusmanovich \cite{ZusSc}, Berman and Kryliouk \cite{BerKry}
the result of Kassel, suitably modified, 
has been extended to
semisimple $\ga$ and  perfect $\ga$ respectively.
But please note  the fact that none of these more general
results on universal central extensions will be used here
to classify the local cocycles.

As a remark aside, not directly related to the topic of
this article, let me add that universal central 
extensions in the 
case  that $\A$ is the algebra of $C^\infty$-functions
on a manifold are studied  for example by Pressley and Segal \cite{PS},
Etingof and Frenkel \cite{EtFren},
Frenkel and Khesin  \cite{FrenKhe}, and Roger \cite{Rogext}.
If $\A$ is a toroidal Lie algebra, i.e. an 
algebra of Laurent polynomials in $N$ variables, the 
universal central extensions and representations of it 
have been studied in more detail.
For example see
Moody, Rao and Yokonoma \cite{MRY90},
Berman and Billig \cite{BB98}, and 
Larsson \cite{Lar00}.

Determining the cohomology of infinite-dimensional
Lie algebras is an important and demanding task.
In the context of classical current algebras let me
mention the work of Feigin \cite{Feicoh}
beside the references already quoted above.
There is a huge literature on the cohomology of vector
field algebras, e.g. see the book of Fuks \cite{Fucohom} which contains
a long list of further references.
For the meromorphic vector field algebras related to higher genus Riemann
surfaces (i.e. the Krichever-Novikov vector field algebras)
the importance of  local  cocycles was realized
by Krichever and Novikov \cite{KNFa,KNFb,KNFc}.
They showed in the case of two marked points for these algebras 
the uniqueness of the local cocycle class.
See \cite{SchlCo} for generalization of the result.
For work in direction of determining the full cohomology
of the Krichever-Novikov vector field algebras see for example
the work of Feigin \cite{Feiicm},  Kawazumi \cite{KawGF}, 
and  Wagemann \cite{Wagkn,Wagdens}.

The structure of this article is as follows.
In \refS{kn} the basic definitions of the
multi-point  higher genus function, vector field, and differential
operator algebras of Krichever-Novikov type are recalled.
The classification results (uniqueness, etc.) of local cocycles
for these algebras are quoted  from \cite{SchlCo}.
In \refS{affine} the definition of the 
higher genus affine algebras 
for general Lie algebras are given.
The above mentioned 
classification results for the reductive case are shown.
In \refS{diff} the algebras  $\D_\ga$  of differential operators
of degree $\le 1$ for the current algebra $\gb$ are studied.
Finally, in \refS{special} the algebras associated to 
$\gl(n)$ and $\sln(n)$ are 
considered as examples .
 
It is a pleasure for me to thank Oleg Sheinman whose
questions inspired me to attack systematically the problem of
almost-graded  central extensions of current algebras, not only for
the current algebras of simple Lie algebras.
He also pointed out to me the importance of the algebras 
$\widehat{\D_{\ga}}$ in
the context of fermionic representations 
\cite{Shferm,Shsc}.
Discussions on universal central extensions 
with Karl-Hermann Neeb and Friedrich 
Wagemann at the {\it Workshop on Lie Theory 2002} in Darmstadt
are acknowledged. Even if the
topic  of the discussion finally did not 
show up in this work directly, it helped to clear up 
my thoughts.
%\newpage
%%%%%%%%%%%%%%%%%%%%%%%%%%%%%%%%%%%
\section{The multi-point algebras of Krichever-Novikov type}\label{S:kn}
%\input  kn.tex
%%%%%%%%%%%%%%%%%%%%%%%%%
%   definition of KN  set-tup
%%%%%%%%%%%%%%%%%%%%%%%%%%%
%    15.7.02// 17.10.
%%%%%%%%%%%%%%%%%%%%%%%%%%%%%%%%%%%%
%\section{The multi-point Krichever-Novikov algebras}\label{S:kn}
%\input kn.tex
%\section{The multi-point Krichever-Novikov algebras}\label{S:kn}
In this section I recall the basic set-up and definitions
for the associative algebra of functions, the
Lie algebra of vector fields, and the Lie algebra of differential
operators of Krichever - Novikov type,
including their central extensions.
In particular, the uniqueness results 
obtained in \cite{SchlCo} and needed further down are
recalled.
The introduction of the current algebra,
of the affine  algebra,
and of the Lie algebra valued differential operator
algebra associated to the current algebra 
will be postponed to the next sections.

\subsection{Geometric set-up and the involved algebras}
%$ $
%\medskip

Let $M$ be a compact Riemann surface of genus $g$, or in terms
of algebraic geometry, a smooth projective curve over $\C$.
Let $N,K\in\N$ with $N\ge 2$ and $1\le K<N$. Fix
$\ 
I=(P_1,\ldots,P_K),$\ {and}\ $O=(Q_1,\ldots,Q_{N-K})$\ 
disjoint  ordered tuples of  distinct points (``marked points'',
``punctures'') on the Riemann surface. 
 In particular, we assume $P_i\ne Q_j$ for every
pair $(i,j)$. The points in $I$ are
called the {\it in-points}, the points in $O$ the {\it out-points}.
Sometimes I will consider $I$ and $O$ simply as sets and use
$A=I\cup O$.

Let $\K$ be the canonical line bundle of $M$.
Its associated sheaf of local sections is the sheaf of
holomorphic differentials.
Following the common practice, I will usually not
distinguish between a line bundle and its associated invertible sheaf
of section.
For every $\l\in\Z$ we consider the bundle
$\ \K^\l:=\K^{\otimes \l}$. Here I use the usual convention:
$\K^0:=\Cal O$ is the trivial bundle, and $\K^{-1}:=\K^*$
is the holomorphic tangent line bundle (resp.
the sheaf  of holomorphic vector fields).
After
fixing a theta characteristics, i.e. a bundle  $S$ with
$S^{\otimes 2}=\K$, we can also consider $\l\in \frac {1}{2}\Z$
with respect to the chosen theta characteristics.
In this article we will only need $\l\in\Z$.
Denote by $\Fl$ the (infinite-dimensional) vector space% 
\footnote{
For $\lambda=\frac 12+\Z$ we should denote the vector space by
$\mathcal{F}^\l_S$
and let $S$ go through all theta characteristics. 
}
of
global meromorphic sections  of $\K^\l$
 which are holomorphic
on $M\setminus A$.

Special cases, which are of particular interest to us, are
the quadratic differentials ($\l=2$),
the  differentials ($\l=1$),
the functions  ($\l=0$), and
the vector fields ($\l=-1$).
The space of functions will be denoted by $\A$ and the
space of vector fields by $\L$.
By  multiplying  sections with functions
we again obtain sections. In this way
the space $\A$ becomes an associative algebra and the spaces $\Fl$ become
$\A$-modules.

The vector fields in $\L$ operate on $\Fl$ by taking
the Lie derivative $L$.
In local coordinates
\begin{equation}\label{E:Lder}
L_e(g)_|:=(\tilde e(z)\fpz)\ldot (\tilde g(z)\dzl):=
\left( \tilde e(z)\pfz {\tilde g}(z)+\l\, \tilde g(z)
\pfz {\tilde e}(z)\right)\dzl \ .
\end{equation}
Here $e\in \L$ and $g\in \Fl$,
and $\tilde e$ and $\tilde g$ are their local representing functions.

The space $\L$ becomes a Lie algebra with respect to \
the Lie derivative \refE{Lder}
and the spaces $\Fl$ become Lie modules over $\L$.
As usual I write $[e,f]$ for the bracket of the vector fields.

For the Riemann sphere ($g=0$) with quasi-global coordinate $z$
 and $I=(0)$ and $O=(\infty)$, the introduced 
function algebra is the algebra of Laurent polynomials 
$\C[z,z^{-1}]$, and the 
vector field algebra is
the Witt algebra, i.e.  the algebra whose universal central extension
is the Virasoro algebra.
We denote for short this situation as the
{\it classical situation}.
%%%%%%%%%%%%%%%%%%%%%%%%%%%%%%%%%%%%%%%%%%%%%%%%%

The vector field algebra $\L$ 
operates on the algebra $\A$ of functions
as derivations.
Hence it is possible to consider the semi-direct sum
$\Do=\A\rtimes \L$.
This Lie algebra is the algebra of 
meromorphic differential operators of degree
$\le 1$ which are holomorphic on $M\setminus A$.
As vector space $\Do=\A\oplus\L$ and the Lie product is given as
\begin{equation}
[(g,e),(h,f)]:=(e\ldot h-f\ldot g,[e,f]).
\end{equation}
There is the short exact sequence of Lie algebras
\begin{equation}\label{E:fdia}
\begin{CD}
0@>>>\A@>i_1>>\Do@>p_2>>\L@>>>0.
\end{CD}
\end{equation}
Obviously, $\L$ is also a  subalgebra of $\Do$.
The vector spaces $\Fl$ become $\Do$-modules by the canonical 
definition
$(g,e)\ldot v=g\cdot v+e\ldot v,\quad v \in\Fl$.
By universal constructions algebras of differential operators of 
arbitrary degree can be considered 
\cite{SchlDiss, Schlct, Schlwed}.
There is another type of algebra of importance, the current algebra.
It will be defined in \refS{affine}.

%%%%%%%%%%%%%%%%%%%%%%%%%%%%%%%%%%%%%%

Let $\rho$ be a meromorphic differential which is holomorphic on 
$M\setminus A$ with exact pole order $1$ at the points in $A$,  given
positive residues at $I$, given negative residues at $O$
(of course obeying the restriction 
$\sum_{P\in I}\res_P(\rho)+
\sum_{Q\in O}\res_Q(\rho)=0$), and purely imaginary periods.
There exists exactly one such $\rho$ (see \cite[p.116]{SchlRS}).
For $R\in M\setminus A$ a fixed point, the function
$u(P)=\Re\int_R^P\rho$ is a well-defined harmonic function.
The family of level lines
$ C_\tau:=\{P\in M\mid u(P)=\tau\},\ \tau\in\R $,
defines a fibration of $M\setminus A$.
Each $C_\tau$ separates the points in $I$ from the points in $O$.
For $\tau\ll 0$ ($\tau\gg 0$) each level line $C_\tau$ is a disjoint union of
deformed circles $C_i$ around the points $P_i$, $i=\iK$ 
(of deformed circles $C_i^*$  around the points $Q_i$, $i=1,\ldots,N-K$).

For $f\in\Fl$ and $g\in\Fl[\mu]$ we have  $f\otimes g\in\Fl[\l+\mu]$.
In particular, for $\mu=1-\l$ 
the form $f\otimes g$ is
a meromorphic differential.
\begin{definition}\label{D:knpair}
The {\it Krichever-Novikov pairing} ({\it KN pairing}) is the
pairing between $\Fl$ and $\Fl[1-\l]$ given by
\begin{equation}\label{E:knpair}
\begin{gathered}
\Fl\times\Fn[1-\l]\ \to\ \C,
\\
\kndual {f}{g}:=\cint{{C_\tau}}f\otimes g
=\sum_{P\in I}\res_{P}(f\otimes g)=
-\sum_{Q\in O}\res_{Q}(f\otimes g),
\end{gathered}
\end{equation}
where $C_\tau$ is any non-singular level line.
\end{definition}
The last equalities follow from the residue theorem.
Note that in \refE{knpair} the integral does not depend on
the level line chosen.
We will call any such level line or any cycle cohomologous to such
a level line a separating
cycle $C_S$.

%%%%%%%%%%%%%%%%%%%%%%%%%%%%%%%%%%%%%%%%%%%%%%%%%%%%%%%%%%%%%
%%%%%%%%%%%%%%%%%%%%%%%%%%%%%%%%%%%%%%%%%%%%%%%%%%%%%%%%%%%%%%%%%%%%
\subsection{Almost-graded structure}
%\medskip
For infinite dimensional algebras 
and their representation theory a graded structure is usually
of importance to obtain structure results.
A typical example is given by the Witt algebra $W$. $W$ admits a
preferred set of basis elements given by
$\{e_n=z^{n+1}\fpz\mid n\in\Z\}$.
One calculates $[e_n,e_m]=(m-n)e_{n+m}$.
Hence $\deg(e_n):=n$ makes $W$ to a graded Lie algebra.

In our more general context the algebras will almost never be graded.
But it was observed by Krichever and Novikov 
in the two-point case that a weaker
concept, an almost-graded structure (they call it a quasi-graded 
structure), will be enough to develop an
interesting  theory of representations (highest weight representations, 
Verma modules, etc.).
\begin{definition}\label{D:almgrad}
(a) Let $\L$ be an (associative or Lie) algebra admitting a direct
decomposition as vector space $\ \L=\bigoplus_{n\in\Z} \L_n\ $.
The algebra $\L$ is called an {\it almost-graded}
algebra if (1) $\ \dim \L_n<\infty\ $ and (2)
there are constants $R$ and  $S$ with
\begin{equation}\label{E:eaga}
\L_n\cdot \L_m\quad \subseteq \bigoplus_{h=n+m+R}^{n+m+S} \L_h,
\qquad\forall n,m\in\Z\ .
\end{equation}
The elements of $\L_n$ are called {\it homogeneous  elements of degree $n$}.
\newline
(b) Let $\L$ be an almost-graded  (associative or Lie) algebra
and $\Cal M$ an $\L$-module with
decomposition $\ \Cal M=\bigoplus_{n\in\Z} \Cal M_n\ $
as vector space. The module $\Cal M$ is called an {\it almost-graded}
module, if
(1) $\ \dim \Cal M_n<\infty\ $, and
(2) there are constants $R'$ and  $S'$ with
\begin{equation}\label{E:egam}
 \L_m \cdot\Cal M_n\quad \subseteq \bigoplus_{h=n+m+R'}^{n+m+S'} \Cal M_h,
\qquad \forall n,m\in\Z\ .
\end{equation}
The elements of $\Cal M_n$ are called 
{\it homogeneous  elements of degree $n$}.
\end{definition}

For the 2-point situation for 
$M$ a higher genus Riemann surface and  $I=\{P\}$, $O=\{Q\}$
with $P,Q\in M$, Krichever and Novikov \cite{KNFa}
introduced an almost-graded structure of the algebras and the modules
by exhibiting
special bases and defining their elements to be the
homogeneous elements.
In \cite{SchlDiss,Schlce} its multi-point
generalization was given, again  by
exhibiting a special basis.
(See also Sadov \cite{Sad} for some results in  similar directions.)

In more detail, 
for fixed $\l$ 
and for  every $n\in\Z$, and $i=1,\ldots,K$, a certain element
$f_{n,p}^\l\in\Fl$ is exhibited  in \cite{SchlDiss,Schlce}.
The $f_{n,p}^\l$ for $p=1,\ldots,K$ are a basis of a
subspace $\Fln$ and it is shown that
$$
\Fl=\bigoplus_{n\in\Z}\Fln\ .
$$
The subspace  $\Fln$ is called the {\it homogeneous subspace of degree $n$}.
%%%%%%%%%%%%%%%%%%%%%%%%%%%%%%%%%%%%%%%%%%%%%%%%%%%%%%%%%

The basis elements are chosen in such a way that they
fulfill the duality relation
with respect to the KN pairing \refE{knpair}
\begin{equation}\label{E:dual}
\kndual {f_{n,p}^\l} {f_{m,r}^{1-\l}}
=\de_{-n}^{m}\cdot
\de_{p}^{r}\ . 
\end{equation}
This implies that the KN pairing is non-degenerate.

We will need as  additional 
information about the elements  $f_{n,p}^\l$ that
\begin{equation}\label{E:ordfn}
\ord_{P_i}(f_{n,p}^\l)=(n+1-\l)-\d_i^p,\quad i=1,\ldots,K .
\end{equation}
The recipe for choosing the orders at the points 
in $O$ is such that up to a scalar multiplication
there is a unique such element which also fulfills \refE{dual}.
After choosing local coordinates $z_p$ at the points
$P_p$ the scalar can be fixed by requiring
\begin{equation}
{f_{n,p}^\l}_|(z_p)=z_p^{n-\l}(1+O(z_p))\left(dz_p\right)^\l\ .
\end{equation}
We  introduce the following
notation:
\begin{equation}\label{E:conc}
A_{n,p}:=f_{n,p}^0,\quad
e_{n,p}:=f_{n,p}^{-1},\quad
\w^{n,p}:=f_{-n,p}^1,\quad
\Omega^{n,p}:=f_{-n,p}^2 \ .
\end{equation}

A detailed analysis \cite{SchlDiss,Schlce} using \refE{dual} yields
\begin{theorem}\label{T:almgrad}
With respect to the above introduced grading the 
associative algebra
$\A$, and the Lie algebras 
$\L$ and $\Do$ are almost-graded algebras and the modules $\Fl$ are
almost-graded modules over them.
In all cases the lower shifts in the degree of the
result (i.e. the numbers $R,R'$ in
\refE{eaga} and \refE{egam}) are zero.
\end{theorem}
Let us abbreviate the  terms of higher degrees as the one
under consideration with the symbol $\hDT$.
By calculating the exact residues in the case of the lower bound
we obtain
\begin{proposition}\label{P:boundary}
%\begin{equation*}
\begin{alignat*}{2}
A_{n,p}\cdot A_{m,r}&=\de_p^r\cdot A_{n+m,r}+\hDT,\quad
&A_{n,p}\cdot f_{m,r}^\lambda&=\de_p^r\cdot f_{n+m,r}^\lambda+\hDT,\quad
\\
[e_{n,p}, e_{m,r}]&=\de_p^r\cdot(m-n)\cdot e_{n+m,r}+\hDT,\quad
&e_{n,p}\ldot f_{m,r}^\lambda&=\de_p^r\cdot (m+\l n)\cdot f_{n+m,r}^\lambda+\hDT.
\end{alignat*}
%\end{equation*}
\end{proposition}
%%%%%%%%%%%%%%%%%%%%%%%%%%%%%%%%%%%%%%%%%%%%%
%%%%%%%%%%%%%%%%%%%%%%%%%%%%%%%%%%%%%%%%%%%%%%%%%%%%%%%%%%%%
%%%%%%%%%%%%%%%%%%%%%%%%%%%%%%%%%%%%%%%%%%%%%%%%%

\subsection{Central extensions and local cocycles}\label{SS:cocyc}
%${ }$

Let $\V$ be a Lie algebra and $\g$  
a Lie algebra two-cocycle of $\V$
with values in $\C$, i.e.
$\g$ is an antisymmetric bilinear form on $\V$ obeying
\begin{equation}\label{E:cocycle}
\g([f,g],h)+\g([g,h],f)+\g([h,f],g)=0,\quad\forall f,g,h\in\V.
\end{equation}
On $\widehat\V:=\widehat\V_\gamma=\C\oplus \V$ a Lie algebra structure
can be defined
by  (with the notation $\widehat{f}:=(0,f)$ and $t:=(1,0)$)
\begin{equation}
[\widehat{f},\widehat{g}]:=\widehat{[f,g]}+\g(f,g)\cdot t,\quad
[t,\widehat\V]=0.
\end{equation}
The element $t$ is a central element.
The algebra $\widehat\V_\gamma$ is called a central extension of $\V$.
Up to equivalence central extensions are classified by the
elements of  $\H^2(\V,\C)$, the second
Lie algebra cohomology space with values in the trivial module $\C$.
In particular, two cocycles $\g_1,\g_2$ define
equivalent central extensions if and only if there exist a
linear form $\phi$ on $\V$ such that
\begin{equation}
\g_1(f,g)=\g_2(f,g)+\phi([f,g]).
\end{equation}
Given a linear form $\phi$ on $\V$ we denote the cocycle obtained
as coboundary by $\phi$ by the symbol
$\delta\Phi$, i.e.
\begin{equation}
\delta\phi(f,g):=\phi([f,g]).
\end{equation}
\begin{definition}\label{D:local}
Let $\V=\bigoplus_{n\in\Z} \V_n$ be an almost-graded Lie algebra.
A cocycle $\g$ for $\V$ is called {\it local}
(with respect to the almost-grading) if there exist $M_1,M_2\in\Z$ with
\begin{equation}\label{E:local}
\forall n,m\in\Z:\quad \g(\V_n,\V_m)\ne 0\implies
M_2\le n+m\le M_1.
\end{equation}
\end{definition}
By defining $\deg(t):=0$ (or any other integer),
the central extension $\widehat\V_\gamma$ is almost-graded if and only if
it is given by a local cocycle $\g$. In this case we call
$\widehat\V_\gamma$ an {\it almost-graded central extension} 
(sometimes also called 
 {\it local central extension}).

In the following we consider cocycles of geometric
origin for the algebras introduced above.  A thorough treatment for them is given in
\cite{SchlCo}.  The proofs of the following statements
and more details  can be found there.

For the abelian Lie algebra $\A$ any antisymmetric bilinear form
will be a 2-cocycle.
Let $C$ be any (not necessarily connected)
differentiable cycle in $M\setminus A$ then
\begin{equation}\label{E:fung}
\g^{(f)}_C:\A\times \A\to\C,\quad
\gamma_C^{(f)}(g,h):=\cinc{C} gdh
\end{equation}
is antisymmetric, hence a cocycle.
Note that by replacing $C$ by any
homologous (differentiable) cycle
one obtains the same cocycle.
The above cocycle is $\L$-invariant, i.e.
\begin{equation}\label{E:linv}
\g^{(f)}_C(e\ldot g,h)+\g^{(f)}_C(g,e\ldot h)=0,
\quad \forall e\in\L,\ \forall g,h\in\A,
\end{equation}
and it is multiplicative, i.e.
\begin{equation}\label{E:mult}
\g(fg,h)+\g(gh,f)+\g(hf,g)=0,
\quad \forall f, g,h\in\A.
\end{equation}

To generalize the
Virasoro-Gelfand-Fuks cocycle 
$\oint (e'''f-ef''')$ for the Witt algebra to
higher genus  vector field algebras $\L$
one has 
first to choose a projective connection
and 
add a counter term to the integrand
to  obtain
a well-defined 1-differential which can be integrated.
Let $R$ be  a global holomorphic projective connection,
see e.g. \cite{SchlCo} for the definition.
For every cycle $C$ (and every $R$) a cocycle  is given by

\begin{equation}\label{E:vecg}
\gamma^{(v)}_{C,R}(e,f):=\cintl{C}\left(\frac 12(\tilde e'''\tilde f-\tilde e\tilde f''')
-R\cdot(\tilde e'\tilde f-\tilde e\tilde f')\right)dz\ .
\end{equation}
Here $e_|=\tilde e\frac {d}{dz}$ and 
$f_|=\tilde f\frac {d}{dz}$ with local meromorphic functions
$\tilde e$ and $\tilde f$.
A different choice of the projective connection
(even if we allow meromorphic projective connections with
poles  possibly at the points in  $A$)
yields a cohomologous cocycle, hence an equivalent central
extension.

These two types of
cocycles can be extended to cocycles {of} $\Do$ by setting
{them} to
be zero if one of the entries are
from the complementary space.
For the vector field cocycles this is clear,
because it is obtained by pulling back the cocycle to $\Do$ via
\refE{fdia}.
For the extension of the function algebra
cocycles the $\L$-invariance  \refE{linv} is crucial.
But there are other independent types of cocycles which mix
functions with vector fields.
To define them
we  first  have  to fix an affine connection $T$
which is holomorphic outside $A$ and has at most a pole of
order one at the point $Q_{N-K}\in O$. For the definition and existence of 
affine connections, see \cite{SchlDiss},\cite{SchlCo},\cite{Shsc}.
The cocycle is defined by zero on elements of the same
type (vector fields or functions) and by
\begin{equation}\label{E:mixg}
\g_{C,T}^{(m)}(e,g):=-\g_{C,T}^{(m)}(g,e):=
\cinc {C}\left(\tilde e\cdot g''+T\cdot (\tilde e\cdot
 g')\right)dz
\end{equation}
for $e\in\L$ and $g\in\A$.
Again, the cohomology class does not depend on the chosen
affine connection.

Now we consider cocycles obtained by integration over a
separating cycle $C_S$.
Instead of $\g_{C_S}$ we will use the symbol $\g_S$.
Clearly, these cocycles can be expressed via residues at the
points in $I$ or equivalently at the points in  $O$.
\begin{proposition}\label{P:local}
\cite{SchlDiss}
The above cocycles \refE{fung}, \refE{vecg} and \refE{mixg}
if integrated over a separating cycle $C_S$
are local.
More precisely,
there exist constants $T_1,T_2,T_3$ such that
for all $m,n\in\Z$
\begin{equation}
\begin{aligned}
\g_{S}^{(f)}(\A_n,\A_m)\ne 0 &\implies  T_1\le m+n\le 0,
\\
\g_{S,R}^{(v)}(\L_n,\L_m)\ne 0 &\implies  T_2\le m+n\le 0,
\\
\g_{S,T}^{(m)}(\L_n,\A_m)\ne 0 &\implies  T_3\le m+n\le 0.
\end{aligned}
\end{equation}
\end{proposition}
If we replace $R$ or $T$ by other meromorphic connections which  have
poles only
at $A$, the cocycle still  will  be local.

One of the main results of \cite{SchlCo} is
\begin{theorem}\label{T:unique}
(a) Every local cocycle of $\A$ which is $\L$-invariant 
is multiplicative and vice versa.
It is given as a
multiple of the cocycle $\g_{S}^{(f)}$.
The cocycle  $\g_{S}^{(f)}$ is cohomologically non-trivial.
\newline
(b) Every local cocycle of $\L$ is cohomologous to
a scalar multiple of $\g_{S,R}^{(v)}$.
The cocycle   $\g_{S,R}^{(v)}$ defines a non-trivial cohomology class.
For every cohomologically non-trivial local cocycle
a meromorphic projective connection $R'$
which is holomorphic outside of $A$ can
be chosen such that the cocycle is equal to a  scalar multiple of
 $\g_{S,R'}^{(v)}$.
\newline
(c) Every local cocycle $\g$ for $\Do$ is a linear combination
of the above introduced cocycles  $\g_{S}^{(f)}$,  $\g_{S,R}^{(v)}$, and
 $\g_{S,T}^{(m)}$ up to coboundary, i.e.
there exists $r_1,r_2,r_3\in\C$ with 
\begin{equation}
\g=r_1 \g_{S}^{(f)}+
r_2 \g_{S,T}^{(m)}+
r_3 \g_{S,R}^{(v)}+\text{coboundary}.
\end{equation}
The 3 basic cocycles are linearly independent in the cohomology space.
If the scalars $r_2$ and $r_3$
in the linear combination are non-zero, then a
meromorphic projective connection $R'$ and an affine
connections $T'$, both  holomorphic outside  $A$,
can be found such that 
$\g=r_1 \g_{S}^{(f)}+
r_2 \g_{S,T'}^{(m)}+
r_3 \g_{S,R'}^{(v)}$.
\end{theorem}
%%%%%%%%%%%%%%%%%%%%%%%%%%%%%%%%%%%%%%%%
Multiplying the cocycle with a nonzero scalar 
corresponds to rescaling the central element. This is an
isomorphism of the central extension.
\begin{corollary}
With respect to the almost-grading up to equivalence and
rescaling (e.g. up to isomorphy) there exist
\newline
(a) a unique almost-graded nontrivial central extension $\Ah$ of $\A$ which is
defined by an $\L$-invariant cocycle,
\newline
(b) a unique almost-graded  nontrivial central extension $\Lh$ of $\L$.
\newline
These central extensions are given by the above geometric cocycles.
\end{corollary}
%%%%%%%%%%%%%%%%%%%%%%%%%%%%
For $\Do$ one obtains up to equivalence a three-dimensional 
family of one-dimensional central extensions, or equivalently,
a three-dimensional central extension.

\begin{remark}
In the context of semi-infinite wedge representations 
associated to the module 
 $\Fl$ (see \cite{SchlDiss},\cite{Schlct},\cite{Schlwed})
as cocycle the linear combination
\begin{equation}
\g_\l=
-(\g_S^{(f)}+\frac {1-2\l}{2}\g_{S,T_\l}^{(m)}+
2(6\l^2-6\l+1)\g_{S,R_\l}^{(v)}),
\end{equation}
with a suitable meromorphic affine connection $T_\l$ and a projective
connection $R_\l$ without poles outside of $A$ and at most
poles of order one at the points in $I$ for $T_\l$ and poles of order two for
 $R_\l$ appears.
It depends on $\l$ and 
defines a central
extension $\Dh_\lambda$ which has a representation 
 on the space of semi-infinite wedge forms.
The cocycle and the action can be extended to a 
central extension $\Dalh$ of the whole differential operator
algebra.
\end{remark}
%%%%%%%%%%%%%%%%%%%%%%%%%%%%%%%%%%%%%%%%%%%%%%%%%%%%%%%%%%%%%
%%%%%%%%%%%%%%%%%%%%%%%%%%%%%%%%%%%%%%%%%%%%%%%%%%%%%%%%%%%%%%%%%%%%
%%%%%%%%%%%%%%%%%%%%%%%%%%%%%%%%%%
%%%%%%%%%%%%%%%%%%%%%%%%%%%%%%%%%%
\section{Affine Algebras}
\label{S:affine}
%\input affine.tex
%\section{Application:  %Cocycles for the affine algebra}
%\label{S:affine}
%\input affine.tex
%  date: 12.9.02// 17.10
%%%%%%%%%%%%%%%%%%%%%%%%%%%%%%%%%%%%%%%%%%%%%%%%%%
In this section I introduce the higher genus analogue
of the 
current algebras and their central extensions, the 
affine algebras (i.e. the Kac-Moody algebras of untwisted 
affine type).
The current algebras  come with an almost-grading and one of the
main goals of this article is to
determine their almost-graded central extensions.
Almost-graded central extensions are also called
local central extensions.
%%%%%%%%%%%%%%%%%%%%%%%%%%%%%%%%%%%%%%%%%%%%%
\subsection{Higher genus current algebras}
Let $\ga$ be an arbitrary  finite-dimensional Lie algebra. 
The {\it multi-point higher genus current algebra} (or  {\it multi-point
higher genus
loop  algebra}) is defined as
\begin{equation}\label{E:curpr}
\gb:=\ga\otimes \A,\quad
\text{with Lie product}\quad
[x\otimes f,y\otimes g]:=[x,y]\otimes (f\cdot g).
\end{equation}
We introduce a grading in  $\gb$ by defining
\begin{equation}
\deg(x\otimes A_{n,p}):=n,\quad n\in\Z,\ p=1,\ldots,K.
\end{equation}

Here $A_{n,p}$ are the basis elements \refE{conc} of the
algebra $\A$.
This makes  $\gb$ to an almost-graded Lie algebra,
$\gb=\oplus_{n\in\Z}(\gb)_n$, with homogenous subspaces
$(\gb)_m=\ga\otimes \A_m$.
To simplify the notation I will sometimes use
$x(g)$ to denote  the element $x\otimes g$.
%%%%%%%%%%%%%%%%%%%%%%%%%%%%%%%%%%%%%%%%%%%%%%%%%

The algebra $\ga$ can be embedded into the algebra $\gb$ as $\ga\otimes 1$.
More precisely, 
\begin{lemma}\label{L:onedec} $1=\sum_{p=1}^K A_{0,p}$.
%\begin{equation}
%1=\sum_{p=1}^K A_{0,p}.
%\end{equation}
In particular, $\ga$ is a subalgebra of $(\gb)_0$.
\end{lemma}
The lemma follows from the calculation of the  basis coefficients of the
function 1, using the Krichever-Novikov duality
\refE{knpair}. See \cite[Lemma~2.6]{SchlShWZ1} 
for details.

\medskip
A Lie algebra $\ga$ fulfilling $\ga':=[\ga,\ga]=\ga$ is called {\it perfect}.
Further down we will need the following simple observation.
\begin{proposition}\label{P:perfect}
If $\ga$ is a perfect Lie algebra then $\gb$ is also perfect.
\end{proposition}
\begin{proof}
The elements of the type $x\otimes f$ with $x\in\ga$ and $f\in\A$ generate
$\gb$. Take any $x\otimes g$, then by perfectness of $\ga$ there exists
$z,y\in\ga$ with $x=[z,y]$. Hence $x\otimes f=[z\otimes f,y\otimes 1]$.
\end{proof}

%%%%%%%%%%%%%%%%%%%%%%%%%%%%%%%%%%%%%%%%%
\subsection{Central extensions}

Next we consider central extensions.
As explained in \refSS{cocyc} central extensions are given via
Lie algebra 2-cocycles with values in the trivial module.
For short they will just be called cocycles.

Examples of central extensions are obtained as follows.
Let $\a$ be a 
fixed invariant, symmetric bilinear form  for the Lie algebra $\ga$,
i.e. a form obeying $\a([x,y],z)=\a(x,[y,z])$.
We do not require it to be non-degenerate.
If $\ga$ is abelian, every  symmetric bilinear form
will be invariant and will do.
If $\ga$ is a simple Lie algebra any multiple of the 
Cartan-Killing form will do.
Choose a  multiplicative cocycle $\gamma$
(see \refE{mult} for the definition) for $\A$.
We set 
$\gh_{\a,\g}=\C\oplus\gb$ as vector space and introduce 
a bracket
\begin{equation}\label{E:daff}
[\widehat{x\otimes f},\widehat{y\otimes g}]=
\widehat{[x,y]\otimes (f g)}+\a(x,y)\gamma (f,g)\cdot t,\qquad
[\,t,\gh_{\a,\g}]=0\ ,
\end{equation}
As usual we
set $\widehat{x\otimes f}:=(0,x\otimes f)$
and $t:=(1,0)$.
\begin{proposition}
The vector space
$\gh_{\a,\g}$ with structure \refE{daff} is a 
Lie algebra and a central extension of $\gb$.
\end{proposition}
\begin{proof}
By  general constructions of central extensions 
(see \refSS{cocyc}) it is enough to
show that 
$\psi(x\otimes g, y\otimes h):=\a(x,y)\g(g,h)$ is a 
cocycle for $\gb$.
The antisymmetry is clear.
Now
\begin{equation}
\psi([x\otimes f,y\otimes g],z\otimes h)=
\psi([x,y]\otimes fg,z\otimes h)
=\a([x,y],z)\g(fg,h).
\end{equation}
By cyclically permuting and adding  up the three terms obtained
we get
\begin{equation}
\begin{gathered}
\a([x,y],z)\g(fg,h)+
\a([y,z],x)\g(gh,f)+
\a([z,x],y)\g(hf,g)=
\\
\a([x,y],z)\cdot\left(\g(fg,h)+\g(gh,f)+\g(hf,g)\right)=0
\end{gathered}
\end{equation}
Here I used the invariance of $\a$ and the 
multiplicative property of $\g$.
\end{proof}

 Let $C$ be any nonsingular curve on $M$ then 
$
\gamma_C(f,g)=\cint{C} fdg
$ 
is a multiplicative  cocycle for $\A$
\cite{SchlCo}.
Hence it defines a central extension $\gh_{\a,C}$ of $\gb$.
%%%%%%%%%%%%%%%%%%%%%%%%%%%%%%%%%%%%%%%%%%%%%%%%%%%%%%%%
\subsection{Local cocycles}
%%%%%%%%%%%%%%%%%%%%%%%%%%%%%%%%%%%%%%%%%%%%%%%%%
Of special importance are cocycles which are local
(see \refD{local}).
As explained there, in this case the almost-grading 
of $\gb$ can be extended to the central extension defined via
this cocycle by setting for the central element $t$,  $\deg t=0$.
From \refT{unique}  we immediately get
\begin{proposition}
Assume that $\a\not\equiv 0$, then 
the cocycle $\a(x,y) \g_C(f,g)$ is local if and only if the 
integration cycle $C$ is a separating cycle $C_S$.
\end{proposition}
In this case $\gh_{\a,S}:=\gh_{\a,C_S}$ 
is an almost-graded central extension of $\gb$.

It should be remarked that everything depends on the 
bilinear form $\a$.
If $\ga$ is simple there is up to multiplication with a
scalar only one such form, the Cartan-Killing form.
In this case it is even a non-degenerate form.

The algebras obtained via a cocycle $\g_{C_S}$ (or more 
general via  $\g_{C}$) 
are called the {\it higher genus (multi-point)
affine Lie algebras (or \KN\ algebras of affine type)}.
In the classical situation (and $\ga$ simple) these are nothing else then the
usual  affine Lie algebras (i.e. the untwisted affine Kac-Moody algebras). 
For higher genus and for the two
point situations such algebras were introduced by 
Krichever and Novikov \cite{KNFa,KNFb,KNFc}, and 
studied in more detail by
Sheinman 
\cite{Shea,Shaff}  and by the author
for the multi-point situation \cite{Schlct, Schlhab}.
%%%%%%%%%%%%%%%%%%%%%%%%%%%%%%%%%%%%%%%%%%%%%%
\begin{proposition}\label{P:alinv}
Let $\ga$  be a finite-dimensional Lie algebra which fulfills the
condition $\ga':=[\ga,\ga]\ne 0$.
Let $\g$ be a local cocycle of $\gb$.
Assume that there is an
invariant symmetric bilinear form $\a$ on $\ga$ 
fulfilling $\a(\ga',\ga)\ne 0$,
and
a bilinear form $\g^{(f)}$ on $\A$ such that $\g$ can be written as
\begin{equation}
\g(x\otimes f,y\otimes z)=\a(x,y)\g^{(f)}(f,g),
\end{equation}
then $\g^{(f)}$ is a multiple of the local cocycle
$\g_{S}^{(f)}$
for the function algebra.
In particular this is the case  if $\a$ nondegenerate.
\end{proposition}
\begin{proof}
First,  $\g^{(f)}$ is obviously antisymmetric and hence a cocycle for
$\A$.
We calculate 
\begin{equation}
\g([x\otimes f,y\otimes g],z\otimes h)=\g([x,y]\otimes f\cdot g,z\otimes h)
=\a([x,y],z)\g^{(f)}(f\cdot g,h).
\end{equation}
For  the cocycle condition for the elements
$x\otimes f,y\otimes g$ and $z\otimes h$
we have to permute this cyclically and add up the results.  We obtain
(using the invariance of $\a$)
\begin{equation}
\a([x,y],z)\left(\g^{(f)}(f\cdot g,h)+
\g^{(f)}(g\cdot h,f)+\g^{(f)}(h\cdot f,g\right)=0.
\end{equation}
By the condition $[\ga,\ga]\ne 0$ and by  $\a(\ga',\ga)\ne 0$
it follows that 
we can find $x,y,z\in\ga$ such that $\a([x,y],z)\ne 0$.
This implies that 
$\g^{(f)}$ is a local and multiplicative cocycle.
\refT{unique} yields the claim.
Note that for $\a$ nondegenerate,  $\a(\ga',\ga)\ne 0$ is 
always true. 
\end{proof}
\begin{proposition}
Let $\g$ be a cocycle of $\gb$ which can be written as 
product
\begin{equation}
\g(x\otimes f,y\otimes z)=\a(x,y)\cdot \g^{(f)}(f,g),
\end{equation}
with $\g^{(f)}$ a non-vanishing multiplicative cocycle
for $\A$ and $\a$ a bilinear form on $\ga$ then
$\a$ is a symmetric invariant bilinear form.
\end{proposition}
\begin{proof}
First we apply the multiplicative property for
$\g^{(f)}$ to the triple $(1,1,f)$ and obtain
$\g^{(f)}(1\cdot 1,f)+
\g^{(f)}(1\cdot f,1)+
\g^{(f)}(f\cdot 1,1)=0$.
This implies that $\g^{(f)}(1,f)=0$ for all $f\in\A$.
From the antisymmetry of $\g$ and $\g^{(f)}$ it follows that 
$\a$ is symmetric.
Now take $f,h\in\A$ such that $\g^{(f)}(f,h)\ne 0$.
We obtain
\begin{equation}
\g([x\otimes f,y\otimes 1],z\otimes h)
=\a([x,y],z)\cdot
\g^{(f)}(f\cdot 1,h).
\end{equation}
Adding this cyclically permuted up, yields
\begin{equation}
0=\a([x,y],z)\cdot \g^{(f)}(f,h)+
\a([y,z],x)\cdot \g^{(f)}(h,f)+
\a([z,x],y)\cdot \g^{(f)}(hf,1).
\end{equation}
Hence 
$\left(\a([x,y],z)-\a([y,z],x)\right)\g^{(f)}(f,h)=0$.
From $\g^{(f)}(f,h)\ne 0$ and the symmetry of $\a$ the invariance
$\a([x,y],z)=\a(x,[y,z])$ follows.
\end{proof}

A cohomology class containing a representing cocycle which is local 
is called  a {\it local cohomology class}.
Note that not all elements in a local cohomology class will be 
local.
But the sum of two local cocycles will be local again.
The subspace of $\H^2(\gb,\C)$ consisting of  local
cohomology classes is denoted by 
$\H^2_{loc}(\gb,\C)$.

%%%%%%%%%%%%%%%%%%%%%%%%%%%%%%%%%%%%%%%%%%%%%%%%%%
\subsection{$\L$-invariant cocycles}
%%%%%%%%%%%%%%%%%%%%%%%%%%%%%%%%%%%%%%%%%%%%%%%
The following definition generalizes the definition of
$\L$-invariance for cocycles of $\A$ to cocycles of
arbitrary current algebras $\gb$.
\begin{definition}\label{D:gblinv}
A cocycle $\g$ of $\gb$ is called
{\it $\L$-invariant} if
\begin{equation}\label{E:glinv}
\g(x(e\ldot g),y(h))+\g(x(g),y(e\ldot h))=0,\qquad
\forall x,y\in\ga,\quad e\in\L,\ g,h\in\A.
\end{equation}
\end{definition}
\begin{proposition}\label{P:cobl}
Let $\psi=\delta\phi$ be a coboundary of $\gb$ which is local and
$\L$-invariant then $\psi=0$.
\end{proposition}
\begin{proof}
The coboundary $\psi$ is given as
$ \psi(x(g),y(h))=\phi([x,y](gh))$ with a linear form
$\phi:\gb\to\C$.
By the $\L$-invariance
\begin{equation*}
0=\psi(x(e\ldot g),y(h))+
\psi(x(g),y(e\ldot h))=
\phi([x,y](e\ldot (gh)).
\end{equation*}
We set $h=1$ and obtain
\begin{equation}\label{E:vanish}
\phi([x,y](e\ldot g))=0,\quad\forall x,y\in\ga,\quad e\in\L,\ g\in\A.
\end{equation}
Assume that $\psi\ne 0$ then there must exist $x,y\in\ga$ and $f\in\A$
such that $\phi([x,y](f))\ne 0$.
Hence, there exists also a basis element $A_{k,r}$ such that
 $\phi([x,y](A_{k.r}))\ne0$.
Assume first that $k>0$ then using the almost-graded structure
(see \refP{boundary})
for the module $\A$ over $\L$ employed
for the pairs of elements $(e_{0,r},A_{k,r})$
we can write  (with a certain $S'\in\N$ and 
coefficients $a'_{h,s}\in\C$)
\begin{equation}
e_{0,r}\ldot A_{k,r}=kA_{k,r}+\sum_{h=k+1}^{k+S'} \sum_{s=1}^K a'_{h,s}
 A_{h,s}.
\end{equation}
The appearing 
$A_{h,s}$ in the sum can be replaced by $(e_{0,s}\ldot A_{h,s})/h+
\hDT$ until we can write
\begin{equation}
e_{0,r}\ldot A_{k,r}=kA_{k,r}+\sum_{h=k+1}^S  \sum_{s=1}^K a_{h,s}
(e_{0,s}\ldot A_{h,s})+D,
\quad\text{with}\ D\in\bigoplus_{h>S}^\infty \A_h,
\end{equation}
and $S$ chosen such that $\psi((\gb)_n,(\gb)_m)=0$ for $n+m>S$.
By the assumed locality of $\psi$ such an $S$ exists.
Now we calculate
\begin{equation*}
0=\phi([x,y](e_{0,r}\ldot  A_{k,r}))=
k\phi([x,y](A_{k,r}))+\sum_{h=k+1}^S  \sum_{s=1}^K a_{h,s}
\phi([x,y](e_{0,s}\ldot A_{h,s}))+\phi([x,y](D)).
\end{equation*}
The first equation follows from \refE{vanish}.
On the r.h.s. the second term will vanish by \refE{vanish} too.
The 3. term can be written as $\phi([x,y](D))=\psi(x\otimes D,y\otimes 1)$.
Using \refL{onedec} we obtain that the 3. term will vanish too.
Hence $k\cdot \phi([x,y](A_{k,r}))=0$ in contradiction to the
assumption.
\newline
If $k=0$ we take for the first step
$e_{-1,r}\ldot A_{1,r}=1\cdot A_{0,r}+\text{higher terms}$.
For the higher terms we can apply $e_{0,s}$ as above.
For $k<0$ we start as above and use instead of $e_{0,s}\ldot A_{0,s}$
the element $e_{-1,s}\ldot A_{1,s}$ to obtain $A_{0,s}$ in the
corresponding intermediate step.
Altogether this leads to a contradiction.
\end{proof}
\begin{proposition}
Let $\g$ be an $\L$-invariant and local cocycle for the
current algebra $\gb$ then
\begin{equation}
\g(x,y(g))=0,\qquad
\forall x,y\in\ga,\ g\in\A.
\end{equation}
\end{proposition}
\begin{proof}
We use the $\L$-invariance \refE{glinv} for $e\in\L$, $g=1$ and $h\in\A$.
It follows $\g(x,y(e\ldot h))=0$ for all $e\in\L$ and  $h\in\A$.
Using the same kind of arguments as in the proof of \refP{cobl},
from the locality of $\g$ and the almost-graded structure the claim follows.
\end{proof}
We call a cohomology class which has an $\L$-invariant
cocycle as an representing element an {\it $\L$-invariant cohomology class}.
By $\H^2_{loc,\L}(\gb,\C)$ we denote the subspace of 
$\L$-invariant and local cocycle classes.
\begin{proposition}\label{P:lloc}
Each class in $\H^2_{loc,\L}(\gb,\C)$ contains a unique representing 
cocycle which is local and $\L$-invariant.
\end{proposition}
\begin{proof}
Let $\g_1$ and $\g_2$ be two local and  $\L$-invariant cocycles
representing the same class.
The difference $\g_1-\g_2$ is a coboundary which is also
local and  $\L$-invariant.
\refP{cobl} shows that it has to vanish.
\end{proof}
%%%%%%%%%%%%%%%%%%%%%%%%%%%%%%%%%%%%%%%%%%%%%%%%%%%%%%%%%%%%%
%%%%%%%%%%%%%%%%%%%%%%%%%%%%%%%%%%%%%%%%%%%%%%%%%%%%%%%%%%%
\subsection{Current algebras of reductive Lie algebras}\label{SS:gred}
%%%%%%%%%%%%%%%%%%%%%%%%%%%%%%%%%%%%%%%%%%%%%%%%%%%%%%%%%
%%%%%%%%%%%%%%%%%%%%%%%%%%%%%%%%%%%%%%%%%%%%%%%%%%%%%%%%
Let $\ga$ be a finite-dimensional reductive Lie algebra.
Recall  (e.g. \cite[I.6.4]{BLG}) that $\ga$ is reductive 
if and only if $\ga$ is the direct sum of Lie algebras
$\ga_0,\ga_1,\ldots,\ga_M$ with $\ga_0$ abelian and 
$\ga_1,\ldots,\ga_M$ simple Lie algebras, i.e.
\begin{equation}
\ga=\ga_0\oplus\ga_1\oplus\cdots\oplus\ga_M.
\end{equation}
In particular, $\ga_0$ is the center of $\ga$,
for the  derived algebra 
we obtain $\ga'=[\ga,\ga]=\ga_1\oplus\cdots\oplus\ga_M$,
and we have $[\ga_i,\ga_i]=\ga_i$, and 
$[\ga_i,\ga_j]=0 $ for $i\ne j$.

For the current algebra  $\gb=\ga\otimes\A$ we obtain the
direct decomposition
\begin{equation}
\gb=\gb_0\oplus\gb_1\oplus\cdots\oplus\gb_M,
\quad\text{with}\ 
\gb_i=\ga_i\otimes \A.
\end{equation}
Furthermore, note that every simple Lie algebra is perfect. Hence using 
\refP{perfect} we obtain 
\begin{equation}\label{E:bperf}
[\gb_i,\gb_i]=\gb_i,\quad
[\gb_i,\gb_j]=0, \ i\ne j,\quad 
\gb'=[\gb,\gb]=\overline{\ga'}=\gb_1\oplus\cdots\oplus\gb_M.
\end{equation}
Let $\g$ be a cocycle for $\gb$ then by restriction 
$\g_i:=\g_{|\gb_i\times\gb_i}$ defines a 
cocycle for the subalgebra and direct summand $\gb_i$ for $i=0,1,\ldots,M$.
Vice versa, given a cocycle $\g_i$ on $\gb_i$ 
we define
the extension $\tilde\g_i$  to $\gb\times \gb$
as $\g_i$ on
$\gb_i\times\gb_i$ and zero otherwise.
\begin{lemma}\label{L:cycdecomp}
(a)   $\tilde\g_i$ is a cocycle for $\gb$ if and only if 
$\g_i$ is a cocycle for $\gb_i$.
\newline
(b) $\tilde\g_i$ is a coboundary for $\gb$ if and only if 
$\g_i$ is a coboundary for $\gb_i$.
\newline
(c) Let $\g$ be a cocycle for $\gb$ then
\begin{equation}
\g(x_i,y_j)=0,\quad\text{for}\ x_i\in\gb_i,\ y_j\in\gb_j,
\ \text{and } i\ne j.
\end{equation} 
\end{lemma}
\begin{proof}
(a) One direction was already explained above.
Assume that $\g_i$ is a cocycle for $\gb_i$.
Clearly, $\tilde\g_i$ is an antisymmetric bilinear form for $\gb$.
We have to show the cocycle condition.
Let $x=\sum_{j=0}^Mx_j$,
 $y=\sum_{k=0}^My_k$, and 
 $z=\sum_{l=0}^Mz_l$ be three elements of $\gb$ with its
unique decomposition with
$x_m,y_m,z_m\in\gb_m$ for $m=0,\ldots, M$.
Using \refE{bperf}
we obtain 
%\begin{equation}
\begin{multline}
\tilde\g_i([x,y],z)=
\tilde\g_i([\sum_{j=0}^Mx_j,\sum_{k=0}^My_k],\sum_{l=0}^Mz_l)
=
\\
\sum_{j,k,l}\tilde\g_i([x_j,y_k],z_l)=
\sum_{j,l}\tilde\g_i([x_j,y_j],z_l)=
\g_i([x_i,y_i],z_i).
\end{multline}
%\end{equation}
Hence from the cocycle condition of $\g_i$ the cocycle condition
for $\tilde\g_i$ follows.
\newline
(b) Let $\g$ be a coboundary for $\gb$, i.e. 
$\g(x,y)=\delta\phi(x,y)=\phi([x,y])$ with 
$\phi:\gb\to\C$ a linear form.
Set $\phi_i=\phi_{|\gb_i}$.
For $x_i,y_i\in\gb_i$ we obtain 
\begin{equation}
\g_i(x_i,y_i)=\g(x_i,y_i)=\delta \phi(x_i,y_i)=
\phi([x_i,y_i])=
\phi_i([x_i,y_i]).
\end{equation}
Hence $\g_i$ will be a coboundary for $\gb_i$.
If $\g_i$ is a coboundary for $\gb_i$ then 
\begin{equation}
\g_i(x_i,y_i)=\delta\phi_i(x_i,y_i) 
=\phi_i([x_i,y_i])
\end{equation}
with a linear form $\phi_i$ on $\gb_i$.
We use $\tilde \phi_i$ to denote $\phi_i$ extended by zero 
on the other components of $\gb$.
But
\begin{equation}
\tilde\g_i(x,y)=\g_i(x_i,y_i)=\phi_i([x_i,y_i])=
\tilde\phi_i([x,y]).
\end{equation}
Hence the claim.
\newline
(c)
Because $i\ne j$ at least one of the indices is $\ne 0$. Assume $i\ne 0$.
Because $[\gb_i,\gb_i]=\gb_i$ we can write
$x_i=[x_i^{(1)},x_i^{(2)}]$ with $x_i^{(1)},x_i^{(2)}\in\gb$.
The cocycle condition for $x_i^{(1)},x_i^{(2)}$ and $y_j$ is
\begin{equation}
\g([x_i^{(1)},x_i^{(2)}],y_i)+
\g([x_i^{(2)},y_j],x_i^{(1)})+
\g([y_j,x_i^{(1)}],x_i^{(2)})=0.
\end{equation}
The last two summands are zero because the commutators are zero.
Hence the first summand which equals $\g(x_i,y_j)$ is also zero.
\end{proof}
\begin{proposition}\label{P:reddec}
(a) Let $\g$ be a
cocycle for $\gb$, $\g_i=\g_{|\gb_i\times\gb_i}$ its restriction
to the $i$-th summand, and $\tilde\g_i$ the extension of $\g_i$  by zero 
to $\gb$ again as above, then
\begin{equation}
\g=\sum_{i=0}^M\tilde\g_i=:\oplus_{i=0}^M\g_i.
\end{equation}
(b) The  above decomposition of $\g$ induces 
\begin{equation}
\mathrm{Z}^2(\gb,\C)=\bigoplus_{i=0}^M
\mathrm{Z}^2(\gb_i,\C),\quad
\mathrm{B}^2(\gb,\C)=\bigoplus_{i=0}^M
\mathrm{B}^2(\gb_i,\C),\quad
\mathrm{H}^2(\gb,\C)=\bigoplus_{i=0}^M
\mathrm{H}^2(\gb_i,\C),
\end{equation}
where $\mathrm{Z}^2$ denotes the vector space of 2-cocycles
and  $\mathrm{B}^2$ denotes the vector space of 2-coboundaries.
\newline
(c)
The cocycle $\g$ is a local (resp. an $\L$-invariant) cocycle
if and only if all $\g_i$, $i=0,1,\ldots,M$ are 
local (resp. $\L$-invariant).
In particular,
\begin{equation}
\H^2_{loc}(\gb,\C)=\bigoplus_{i=0}^M \H^2_{loc}(\gb_i,\C),\qquad
\H^2_{loc,\L}(\gb,\C)=\bigoplus_{i=0}^M \H^2_{loc,\L}(\gb_i,\C).
\end{equation}
\end{proposition}
\begin{proof}
Per construction $\g$ coincides with $\sum_i\tilde\g_i$ on
$\gb_j\times\gb_j$ for all $j$.
By \refL{cycdecomp} both also coincide on  
$\gb_j\times\gb_k$ for all $j,k$ with $j\ne k$.
Hence (a).
Part (b) is a reformulation of \refL{cycdecomp} (a) and (b).
Clearly locality and $\L$-invariance is a condition which has to 
be checked on each component $\gb_i$ 
separately. Hence (c) follows.
\end{proof}
One way to obtain central extensions of $\gb$ was to choose a 
symmetric invariant bilinear form $\a$ for $\ga$ and use the 
defining equation \refE{daff}.
With the same type of argument as above any such $\a$ for reductive
$\ga$ can be decomposed as
\begin{equation}
\a=\oplus_{i=0}^M\a_{i},\quad
\text{with}\ 
\a_i=\a_{|\ga_i\times\ga_i}
\end{equation}
symmetric invariant bilinear forms on $\ga_i$.
Vice versa, every such sum of $\a_i$ gives a 
symmetric invariant bilinear form $\a$ on $\ga$.
Note also that $\a(\ga_i,\ga_j)=0$ for $i\ne j$.
Recall that for a simple Lie algebra $\ga_i$ over $\C$ there is
up to multiplication with a scalar $r_i$ only one such
invariant form, the Cartan-Killing form $\b_i$ of $\ga_i$.
The form $\b_i$ is non-degenerate.
We obtain $\a=\a_0+\sum_{i=1}^Mr_i\b_i$ with $r_i\in\C$.
For the abelian Lie algebra $\ga_0$ every symmetric bilinear
form is invariant. If $\dim\ga_0=n$ then the space of
symmetric bilinear forms on $\ga_0$ is $\frac {n(n+1)}{2}$-dimensional.
Altogether we obtain that the space of symmetric invariant 
bilinear forms on $\ga$ is 
$\dfrac {n(n+1)}{2}+M$-dimensional.

\medskip
Our aim is to determine all local (and $\L$-invariant) cocycles
for the current algebra $\gb$ of a reductive Lie algebra $\ga$.
To reach this goal we first have to study the 
cocycles for the current algebras
of the simple, semisimple, and abelian summands of $\ga$.

%%%%%%%%%%%%%%%%%%%%%%%%%%%%%%%%%%%%%%%%%%%%%%%%%%%%%%%%%%%%%%
%%%%%%%%%%%%%%%%%%%%%%%%%%%%%%%%%%%%%%%%%%%%%%%%%%%%%%%%%%%%%
\subsection{Cocycles for the simple case}

\begin{theorem}\label{T:simple}
(a) Let $\ga$ be a finite-dimensional simple Lie algebra, then every
local cocycle of the current algebra $\gb=\ga\otimes \A$ is
cohomologous to a cocycle given by
\begin{equation}\label{E:simple}
\gamma(x\otimes f,y\otimes g)=r\cdot \frac{\b(x,y)}{2\pi\i}\int_{C_S}fdg,\quad
\text{with}\ r\in\C,
\end{equation}
and $\b$ the Cartan-Killing form of $\ga$.
In particular,
$\g$ is cohomologous to a local and  $\L$-invariant cocycle.
\newline
(b) If the cocycle is 
already local and  $\L$-invariant, 
then it coincides with the cocycle
\refE{simple} with  $r\in\C$ suitable chosen.
\end{theorem}
\begin{proof}
Kassel \cite{Kaskd} proved that the algebra $\gb=\ga\otimes\A$ for any
commutative algebra $\A$ over $\C$ and any  simple 
Lie algebra $\ga$ admits
a universal central extension. It is given by
\begin{equation}
\gh^{\, univ}=(\Omega^1_{\A}/d\A)\oplus\gb
\end{equation}
with Lie structure
\begin{equation}
[x\otimes f, y\otimes g]=[x,y]\otimes fg+\b(x,y)\overline{fdg},\qquad
[\Omega^1_{\A}/d\A,\gh^{\, univ}]=0.
\end{equation}
Here $\Omega^1_{\A}/d\A$ denotes the vector space of K\"ahler differentials
of the algebra $\A$, and $\b$ is the Cartan-Killing form.
The elements in $\Omega^1_{\A}$ can be given as $fdg$ with $f,g\in\A$, 
and $\overline{fdg}$ denotes its class modulo $d\A$.
This universal extension is not necessarily one-dimensional.
Let $\gh$ be any one-dimensional central extension of $\gb$. It will  be given
as a quotient of $\gh^{\,univ}$. Up to equivalence it can be given by
a Lie homomorphism $\Phi$
\begin{equation}
\begin{CD}
\gh^{\,univ}=\Omega^1_{\A}/d\A\oplus\gb
@>\ \Phi=(\varphi,id)\ >>\gh=\C\oplus\gb
\end{CD}
\end{equation}
with a linear form $\varphi$ on $\Omega^1_{\A}/d\A$.
The structure of $\gh$ is then equal to
\begin{equation}
[x\otimes f, y\otimes g]=[x,y]\otimes fg+\b(x,y)\varphi(\overline{fdg}) t,\qquad
[t,\gh]=0.
\end{equation}
In our situation $M\setminus A$ is an affine curve and 
$\Omega^1_{\A}/d\A$ is the first cohomology group 
of the complex of
meromorphic functions on $M$ which are holomorphic on $M\setminus A$
(similar arguments can be found in an article by Bremner \cite{Bremce}).
By Grothendieck's algebraic deRham theorem \cite[p.453]{GH}
the cohomology of the complex is isomorphic to 
 the singular cohomology of $M\setminus A$.
Such a linear form $\varphi$ can be given by choosing a linear
combination of cycle classes in $M\setminus A$ and integrating
the differential class $\overline{fdg}$ over this combination.
By \refT{unique} the locality implies that the combination
is a multiple of the separating cocycle.
This shows 
that the given cocycle is cohomologous to \refE{simple}. But
this cocycle is local and $\L$-invariant. Hence, (a).
By \refP{lloc} it is the only  local and $\L$-invariant cocycle in his
class.
Hence, (b).
\end{proof}

\begin{corollary}
For $\ga$ simple,
up to equivalence of extensions and rescaling of the central element there is
a unique non-trivial almost-graded  central extension
$\gh$, its higher genus
multi-point current algebra $\gb$.
It is given by the cocycle \refE{simple}.
\end{corollary}
\begin{corollary}
For $\ga$ simple
\begin{equation}
\H^2_{loc}(\gb,\C)=
\H^2_{loc,\L}(\gb,\C).
\end{equation}
This space is one-dimensional.
\end{corollary}
%%%%%%%%%%%%%%%%%%%%%%%%%%%%%%%%%%%%%%%%%%%%%%%%%%%%%%%%%%%%
%%%%%%%%%%%%%%%%%%%%%%%%%%%%%%%%%%%%%%%%%%%%%%%%%%%%%%%%%%%
\subsection{Cocycles for the semisimple case}
%%%%%%%%%%%%%%%%%%%%%%%%%%%%%%%%%%%%%%%%%%%%%%%%%%%%%%%%%%
Let $\ga$  be a semisimple Lie algebra. It can 
be written as direct sum of
simple Lie algebras $\ga_i$, i.e. $\ga=\ga_1\oplus\cdots\oplus \ga_M$.
\begin{theorem}\label{T:ss}
(a) For every local cocycle $\g$ for the current algebra $\gb$ of
a semisimple Lie algebra $\ga$ there exists a symmetric 
invariant bilinear form $\a$ for $\ga$ such that $\g$ is cohomologous to
\begin{equation}\label{E:ssimple}
 \gamma'_{\a,S}(x\otimes f,y\otimes g)=\frac{\a(x,y)}{2\pi\i}\int_{C_S}fdg.
\end{equation}
In particular,
$\g$ is cohomologous to a local and  $\L$-invariant cocycle.
\newline
(b) If the cocycle is 
already local and  $\L$-invariant, 
then it coincides with the cocycle \refE{ssimple}.
\newline
(c)
\begin{equation}
\dim\H^2_{loc}(\gb,\C)=\dim\H^2_{loc,\L}(\gb,\C)=M,
\end{equation}
where $M$ is the number of simple summands  of $\ga$.
\end{theorem} 
\begin{proof}
As explained in \refSS{gred} every local 
cocycle for $\gb$ can be decomposed 
as $\g=\oplus_{i=1}^M\g_i$ with 
$\g_i$ local cocycles on the simple summands $\g_i$.
Denote by $\b_i$ the Cartan-Killing form of $\ga_i$. 
By \refT{simple} the $\g_i$ is cohomologous to
\begin{equation}
\gamma'_i(x_i\otimes f,y_i\otimes g)=
r_i\cdot \frac{\b_i(x_i,y_i)}{2\pi\i}\int_{C_S}fdg,\quad
\text{with}\ r_i\in\C.
\end{equation}
We set  $\a:=\sum_{i=1}^Mr_i\b_i$ which is defined 
for $x=\sum_j x_j$, $y=\sum_k y_k$ by
$a(x,y):=\sum_{i=1}^M r_i\b_i(x_i,y_i)$.
The form $\a$ is a symmetric invariant bilinear form 
and
\begin{equation}
\gamma_{\a,S}(x\otimes f,y\otimes g)=\frac{\a(x,y)}{2\pi\i}\int_{C_S}fdg
\end{equation}
is a cocycle which is cohomologous to $\g$
(use \refP{reddec}).
This shows (a). Part (b) follows from 
\refP{lloc}.
All linear combinations of the Cartan-Killing forms $\b_i$ give the
whole (local) cohomology space. Hence (c). This follows also from
\refP{reddec} (b).
\end{proof}
\begin{corollary}
The space of equivalence classes of almost-graded central extensions
of the current algebra $\gb$ of a semisimple Lie algebra $\ga$ is in 
1:1 correspondence with the space of symmetric invariant bilinear forms
for $\ga$. Its dimension is the number of simple summands of $\ga$.
\end{corollary}
%%%%%%%%%%%%%%%%%%%%%%%%%%%%%%%%%%%
\begin{remark}
In general, the claim of the above theorem is not true for
$\ga$ not semisimple. 
As a nontrivial example take $\ga=\gl(n)$. Let $\a$ be any
nonvanishing symmetric invariant bilinear form for $\gl(n)$.
For any antisymmetric bilinear form $\psi$ on $\A$
the form 
\begin{equation*}
\gamma(x\otimes f,y\otimes g):=\a(x,y)\psi(f,g)
\end{equation*}
defines a cocycle for $\overline{\gl}(n)$. 
But $\psi$ can be chosen to be local without being
a geometric cocycle.
We have to require the cocycle to be 
$\L$-invariant, to obtain a statement corresponding to the above
theorem
at least in the reductive case, see \refSS{cred}.
\end{remark}
%%%%%%%%%%%%%%%%%%%%%%%%%%%%%%%%%%%%%%%%%%55

%%%%%%%%%%%%%%%%%%%%%%%%%%%%%%%%%%%%%%%%%%%%%%%%%%%%%%%%%%%%
%%%%%%%%%%%%%%%%%%%%%%%%%%%%%%%%%%%%%%%%%%%%%%%%%%%%%%%%%%%
\subsection{Cocycles for the abelian case}
%%%%%%%%%%%%%%%%%%%%%%%%%%%%%
Let $\ga$ be a finite-dimensional abelian Lie algebra.
Note  that in the abelian case there exists no non-trivial
coboundaries. Hence two different cocycles will never be
cohomologous. Already in the one-dimensional abelian case, i.e. 
$\ga=\C$, $\gb=\C\otimes \A=\A$, if one wants to obtain
uniqueness results, it is necessary to 
allow 
only $\L$-invariant cocycles.
If $\dim\ga\ge 2$ an additional phenomena will show up.
In contrast to the semisimple case the Lie algebra $\ga$ will have
non-trivial central extensions.
Those are given by  
alternating 
bilinear forms on $\ga$.
If we would allow arbitrary central extensions of $\gb$ then also
central extensions induced by extensions of $\ga$ would show up.
It will turn out, that the required $\L$-invariance will 
exclude cocycles coming from $\ga$.
\begin{lemma}\label{L:split}
Let $\g$ be a local and $\L$-invariant cocycle for $\gb$.
\newline
Set $\gxy(f,g):=\g(x\otimes f,y\otimes g)$ for fixed $x,y\in\ga$.
Then there exists $\a_{xy}\in\C$ with
\begin{equation}\label{E:split}
\gxy(f,g)=\a_{xy}\g_S^{(f)}(f,g).
\end{equation}
\end{lemma}
\begin{proof}
For fixed $x,y$ the form $\gxy$ will be a bilinear form on $\A$.
The claim follows with the same kind of arguments as in \cite{SchlCo}
used for the $\L$-invariant local cocycles for the function algebra $\A$
to show that they are multiples of $\g_S^{(f)}$.
Any reference in the proof there to the antisymmetry of the
bilinear form giving the 
cocycle can be avoided.
Hence,  $\gxy$ has to be a 
multiple of the function algebra cocycle 
$\g_S^{(f)}$.
In particular the antisymmetry will follow.
For the convenience of the reader the details are given in an
appendix to  this article.
\end{proof}
\begin{theorem}\label{T:abel}
Let $\ga$ be an abelian Lie algebra of dimension $n$. Then
\begin{equation}
\dim\H^2_{loc,\L}(\gb,\C)=\frac {n(n+1)}{2},
\end{equation}
and the cocycles are given by
\begin{equation}\label{E:acent}
\gamma_\a(x\otimes f,y\otimes g)=\frac{\a(x,y)}{2\pi\i}\int_{C_S}fdg,
\end{equation}
where $\a$ is an arbitrary symmetric bilinear form for $\ga$.
\end{theorem}
\begin{proof}
For $n=1$ this is the result for $\gb=\A$ obtained in 
\cite{SchlCo} and here recalled in \refT{unique}.
Clearly, given any such $\a$ we obtain via \refE{acent} a 
cocycle which is local and $\L$-invariant.
Now let $\g$ be a local and $\L$-invariant cocycle.
From \refL{split} we know that we can write
$\g(x\otimes f,y\otimes g)=\a_{xy}\cdot\g_S^{(f)}$.
The map $\alpha:(x,y)\to \a_{x,y}$ is bilinear. By the antisymmetry of 
$\g$, and $\g_S^{(f)}$, and the non-vanishing of the latter,
$\a$ has to be a symmetric form on $\ga$.
Hence, $\g$ is  indeed of the form \refE{acent}.
But there are exactly $\frac {n(n+1)}{2}$ 
symmetric bilinear forms which are linearly independent.
The corresponding $\g_{\a}$s will stay linearly independent.
This shows the formula for the dimension.
\end{proof}
%%%%%%%%%%%%%%%%%%%%%%%%%%%%%%%%%%%%%%%%%%%%%%%%%%%%%%%%%%%%
%%%%%%%%%%%%%%%%%%%%%%%%%%%%%%%%%%%%%%%%%%%%%%%%%%%%%%%%%%%
\subsection{Cocycles for the reductive  case}\label{SS:cred}
%%%%%%%%%%%%%%%%%%%%%%%%%%%%%%%%%%%%%%%%%%%%%%%%%%%%%%%%%%%
Let $\ga=\ga_0\oplus\ga_1\oplus\cdots\oplus\ga_M$ 
the decomposition of the reductive Lie algebra $\ga$ into its
abelian and simple summands.
\begin{theorem}\label{T:red}
(a) Given a cocycle $\g$ for $\gb$ which is local, and 
whose restriction $\g_0$ to the abelian summand 
$\gb_0$ is $\L$-invariant, then there
exists a symmetric invariant bilinear form $\a$ for $\ga$ such that
$\g$ is cohomologous to 
\begin{equation}
 \gamma'_\a(x\otimes f,y\otimes g)=\frac{\a(x,y)}{2\pi\i}\int_{C_S}fdg.
\end{equation}
Vice versa, every such $\a$ determines a local and $\L$-invariant
cocycle.
\newline
(b) If the cocycle $\g$ is 
already local and  $\L$-invariant for the whole $\gb$, 
then it coincides with the cocycle $\g'_a$.
\newline
(c)
\begin{equation}
\dim\H^2_{loc,\L}(\gb,\C)=\frac {n(n+1)}{2}+M.
\end{equation}
\end{theorem}
\begin{proof}
Again (a) is obtained by restricting the cocycle to the
semisimple and the abelian summand.
Using \refT{ss} and \refT{abel} 
together with \refP{reddec} we obtain the results.
\end{proof}
\begin{corollary}
Those equivalence classes of 
almost-graded central extensions of the current algebra
$\gb$ of a reductive Lie algebra $\ga$,
whose corresponding cocycles can be given by $\L$-invariant cocycles, 
correspond 
1:1 to the space of symmetric invariant bilinear forms $\a$  
on $\ga$.
This space has the dimension
$\dfrac {n(n+1)}{2}+M$.
\end{corollary}

%%%%%%%%%%%%%%%%%%%%%%%%%%%%%%%%%%%%%%%%%%%%%%%%%%%%%%%%%%%%
%%%%%%%%%%%%%%%%%%%%%%%%%%%%%%%%%%%%%%%%%%%%%%%%%%%%%%%%%%%
%%%%%%%%%%%%%%%%%%%%%%%%%%%%%%%%%%%%%%
\section{Algebras of differential operators}
\label{S:diff}
%\input diff.tex
%\section{Algebras of differential operators}
%\label{S:diff}     17.10.02
%\input diff.tex
%%%%%%%%%%%%%%%%%%%%%%%%%%%%%%%%%%%
In \refS{kn} the algebra of differential operators 
of degree $\le 1$ associated to  the
function algebra $\A$  was introduced as semi-direct sum of 
 $\A$ with the vector field algebra $\L$.
In this section the construction will be extended to the case when
$\A$ is replaced by 
the current algebra $\gb$ of 
a general finite-dimensional  Lie algebra 
$\ga$.
Again, central extensions are studied
and classification results for local cocycles are given.
The algebras play an important role in the context of
Casimir operators in 
fermionic representations, (e.g. see \cite{Shferm,Shsc}).
The fermionic representations generalize the semi-infinite 
wedge representations which were studied in 
\cite{KNFb},\cite{SchlDiss},\cite{Schlwed}.
%%%%%%%%%%%%%%%%%%%%%%%%%%%%%%%%%%%%%%%%%%%%%%

\subsection{Differential operator algebras}

Let $\ga$ be an arbitrary finite-dimensional Lie algebra, $\gb$ the
associated current algebra.
Set $\D_\ga:=\gb\oplus\L$ as vector space.
By 
\begin{equation}
e\ldot (x\otimes g):=x\otimes(e\ldot g),\qquad
e\in\L,\ x\in\ga,\  g\in\A,
\end{equation}
an  $\L$-module structure is defined on the
space $\gb$.
\begin{definition}
The Lie algebra $\D_\ga$ is the semi-direct 
sum structure on $\gb\oplus\L$ induced by
$[e,x\otimes g]:=x\otimes(e\ldot g)$.
\end{definition}
In more detail, the product is given as
(for $e,f\in\L$, $x,y\in\ga$, $g,h\in\A$)
\begin{equation}\label{E:dstruct}
[(x\otimes g,e),(y\otimes h,f)]=
(\;[x,y]\otimes(g\cdot h+g\cdot(e\ldot h)-(f\ldot g)\cdot h)\;,\; [e,f]\;).
\end{equation}
There is the short exact sequence of Lie algebras
\begin{equation}\label{E:dsequ}
\begin{CD}
0@>>>\gb@>i_1>>\D_\ga@>p_2>>\L@>>>0.
\end{CD}
\end{equation}

Using the almost-gradings of $\gb$ and of $\L$ we obtain
an almost-grading for $\D_\ga$ by taking 
$(\D_\ga)_m:=(\gb)_m\oplus \L_m$ for $m\in\Z$ 
as homogenous subspaces.
From \refE{dstruct} one immediately verifies that 
this defines an  almost-graded structure.
For $\ga=\C$ (and hence $\gb=\A$) we recover the algebra
$\D$ as special case.
%%%%%%%%%%%%%%%%%%%%%%%%%%%%%%%%%%%%%%%

\subsection{Cocycle conditions}

%%%%%%%%%%%%%%%%%%%%%%%%%
Next we consider cocycles and central extensions of $\D_\ga$.
By bilinearity the cocycle condition \refE{cocycle} can be equivalently
formulated by considering cocycle conditions for triples of
elements of ``pure types'',
i.e. elements which are either vector fields or elements from $\gb$.
\begin{enumerate}
\item
If all three are vector fields, the condition is that the
cocycle defines by restriction  a cocycle for $\L$.
\item
If all three elements are from $\gb$, the condition is
that the cocycle defines by restriction a cocycle for $\gb$.
\item
Now let $e,f\in\L$ and $x(g)\in\gb$, then using the Lie structure
we obtain
\begin{equation}\label{E:cefA}
\g([e,f],x(g))-\g(e,x(f\ldot g))+\g(f,x(e \ldot g))=0.
\end{equation}
\item
If  $e\in\L$ and $x(g),y(h)\in\gb$ then
\begin{equation}\label{E:ceAB}
\g(x(e\ldot g),y(h))-\g(e,[x,y](gh))+\g(x(g),y(e\ldot h))=0.
\end{equation}
\end{enumerate}
An  antisymmetric form on $\gb$ will be a 
cocycle if and only if all 4 conditions are fulfilled.

\subsection{Extensions of cocycles}\label{SS:ext}

In the following 
for a cocycle of a subalgebra $V$ of $\Do_\ga$
I will use the expression ``extended by zero on the complementary space''.
By this I mean that 
the extended bilinear form coincides with the 
cocycle on pairs of elements from $V$  
and will be set to zero if any of the two entries 
in the bilinear form is from 
the complementary subspace. 
The rest follows from bilinear extension.
Note that this extension is not necessarily  a cocycle
for $\D_\ga$. 
%%%%%%%%%%%%%%%%%%%%%%%%%%%%%%%%%%%%%%%%%%%%%%%%%
\begin{proposition}
$ $
\newline
(a) Every cocycle of $\L$ can be extended by zero on $\gb$ to 
a cocycle of $\D_\ga$.
\newline
(b)
A cocycle of $\gb$ can  be extended by zero on $\L$ 
to a cocycle of $\D_\ga$ if and only
if the cocycle is $\L$-invariant (see \refD{gblinv}).
\newline
(c) Let $\a$ be an invariant symmetric bilinear form for $\ga$.
The geometric cocycles 
\begin{equation}
\g_{\a,C}(x(g),y(h))=\a(x,y)\cdot \g_{C}^{(f)}(g,h)
\end{equation}
of $\gb$ can be extended by zero to a cocycle for $\D_\ga$.
\end{proposition}
\begin{proof}
From the above 4 separate cocycle 
conditions only the first one is of relevance
for a cocycle of $\L$ which is extended by zero.
This is exactly the condition that it is a cocycle for $\L$.
This follows  also  from the
sequence \refE{dsequ} because the extended cocycle is nothing but the
pull-back by $p_2$. Hence, (a).
The situation is  different for cocycles of $\gb$.
If we set the cocycle to be zero outside of $\gb$, then
(1), (2), and (3) are automatically fulfilled.
From \refE{ceAB} it follows that
it will define a cocycle if and only if
\begin{equation}
\g(x(e\ldot g),y(h))+\g(x(g),y(e\ldot h))=0.
\end{equation}
But this is  \refD{gblinv} of
$\L$-invariance for cocycles of $\gb$. Hence, (b).
Part (c) follows from the fact that 
$\g_{C}^{(f)}$ is $\L$-invariant, hence the same is true for 
$\g_{\a,C}(x(g),y(h))$.
\end{proof}

\subsection{Cocycles of mixing type}\label{SS:mix}

Now we consider  cocycles $\gamma$ of pure mixing type.
By a cocycle of pure mixing type I mean a cocycle which vanishes
on pairs of  elements of the same type.
From \refE{ceAB} it follows that
\begin{equation}\label{E:mcomm}
\gamma(e,[x,y](g))=0,\quad\forall x,y\in\ga,\ e\in\L,\ g\in\A,
\end{equation}
is a necessary condition for $\g$ to be a cocycle of pure mixing type.
\begin{proposition}\label{P:semi}
For a perfect Lie algebra $\ga$  (i.e. $\ga=[\ga,\ga]$)
 there exist no non-vanishing cocycles
of $\Do_{\ga}$ of pure
mixing type.
\end{proposition}
\begin{proof}
By assumption we
 can express every $z\in \ga$ as $z=[x,y]$.
Hence by \refE{mcomm} the claim follows.
\end{proof}
Further down we will need the following result.
\begin{lemma}\label{L:cobmix}
Let $\g$ be a cocycle of pure mixing type which is
a coboundary.
Then there exists a linear
form  $\psi:\gb\oplus\L\to\C$ which vanishes on $\L$
with  $\g=\delta\psi$.
\end{lemma}
\begin{proof}
As a coboundary $\g=\delta\psi'$
with $\psi':\gb\oplus\L\to\C$
a linear form. We decompose $\psi'=\psi'_{|\gb}\oplus\psi'_{|\L}$.
For all $e,f\in\L$ we obtain
\begin{equation}
0=\g(e,f)=\delta\psi(e,f)=\psi([e,f]),
\end{equation}
by the requirement of pure mixing type.
Hence, $\psi'_{|[\L,\L]}\equiv 0$.
We extend $\psi'_{|\gb}$ by zero on $\L$ to a linear form $\psi$.
Note that $\psi([e,x(g)])=\psi(x(e\ldot g))=\psi'([e,x(g)])$
and  $\psi([x(g),y(h)])=\psi'([x(g),y(h)])=0$.
Hence we get $\g=\delta \psi$ with a linear form $\psi$ of
the required kind.
\end{proof}
Let $\g$  be a cocycle of pure mixing type
which is of the form
$\g(e,x(g))=\phi(x)\g^{(m)}(e,g)$ with $\phi\in\ga^*$
a linear form on $\ga$, and
$\g^{(m)}$ bilinear on $\L\times\A$.
Condition
\refE{mcomm} implies that
either $\g^{(m)}\equiv 0$ or that $\phi([x,y])=0$ for all $x,y\in\ga$.
By Condition \refE{cefA} for
$x\in\ga$, $e,f\in\L$ and $g\in\A$  the relation
\begin{equation}\label{E:mixl}
\phi(x)\left(\g^{(m)}([e,f],g)-\g^{(m)}(e,f\ldot g)+
\g^{(m)}(f,e\ldot g)\right)=0
\end{equation}
follows.
Excluding the  trivial case $\phi\equiv 0$ this implies that
$\g^{(m)}$ is a mixing cocycle for $\D$.
Vice versa, every bilinear form $\g(e,x(g))=\phi(x)\g^{(m)}(e,g)$
with $\phi\in\ga^*$ such that  $\phi_{|[\ga,\ga]}=0$, and
$\g^{(m)}$ a mixing cocycle of $\Do$ defines a cocycle of
 $\D_\ga$.
%%%%%%%%%%%%%%%%%%%%%%%%%%%%%%%%%%%%%%%%%%%%%%
\subsection{The general case}\label{SS:dgen}
Without assuming further properties of the finite-dimensional
Lie algebra $\ga$ we obtain the following theorem.
\begin{theorem}\label{T:gdiff}
Let $\ga$ be a finite-dimensional  Lie algebra.
\newline
(a) Let $\g^{(f)}$ be any $\L$-invariant cocycle of the
function algebra, $\g^{(v)}$ any cocycle of the
vector field algebra $\L$, and $\g^{(m)}$ any mixing
cocycle of $\D$. Furthermore, let $\a$ be any symmetric invariant
bilinear form on $\ga$ and $\phi$ any linear form on $\ga$ which
vanishes on the derived subalgebra $\ga'=[\ga,\ga]$, then
\begin{multline}\label{E:sum3a}
\g((x(g),e),(y(h),f))=
r_1\a(x,y)\g^{(f)}(g,h)
\\
+r_2\big(\phi(y)
\g^{(m)}(e,h)-\phi(x)\g^{(m)}(f,g)\big)+r_3\g^{(v)}(e,f)
\qquad
\end{multline}
for arbitrary $r_1,r_2,r_3\in\C$ defines a cocycle of $\D_\ga$.
\newline
(b)
Let $\g$ be a local cocycle of $\D_\ga$.
Assume that it can be written as
\begin{multline}\label{E:sum3b}
\g((x(g),e),(y(h),f))=
r_1\a(x,y)\tilde \g^{(f)}(g,h)
\\
+r_2\big(\phi(y)
(\tilde\g^{(m)}(e,h)-\phi(x)\tilde \g^{(m)}(f,g)\big)+
r_3\tilde\g^{(v)}(e,f)
\qquad
\end{multline}
with a symmetric invariant  bilinear form $\a$
for $\ga$
fulfilling $\a(\ga',\ga)\ne 0$,
$\phi$ a
linear form on $\ga$, $\tilde \g^{(f)}$
a bilinear form on $\A$, $\tilde\g^{(v)}$
a bilinear form on $\L$ and $\tilde \g^{(m)}$ bilinear form on $\L\times\A$,
and $r_1,r_2,r_3\in\C$.
Assume either $[\ga,\ga]\ne 0$ or that  $\tilde \g^{(f)}$ is $\L$-invariant.
Then,  ignoring all terms which are identically zero, the
remaining forms are multiples of the corresponding unique local cocycles
introduced in \refSS{cocyc}, respectively are coboundaries,
 and the  linear form $\phi$ vanishes on
$\ga'=[\ga,\ga]$.
\end{theorem}
\begin{proof}
Part (a) follows from the discussion  in \refSS{mix} and \refSS{ext}.
To see (b) we note that the restriction of $\g$ to $\L$ defines a local
cocycle hence for this part \refT{unique}
shows the correct form.
The restriction of the cocycle to $\gb$ also defines a cocycle.
Hence if  $[\ga,\ga]\ne 0$ we can use \refP{alinv} and obtain the claim
for this part.
In case $[\ga,\ga]=0$ we  assumed
$\L$-invariance of $\tilde\g^{(f)}$ and 
by  locality of the cocycle \refT{unique} shows the claim
directly.
Hence, the first term and the last term define cocycles of $\Do_\ga$. This
implies that the second term is also a cocycle.
In particular, it is a pure mixing cocycle.
For the mixing cocycle the relation \refE{mixl}
shows that either $\phi$
 is identical zero (and the term does not appear at all)
or $\tilde \g^{(m)}$ defines a local cocycle for $\D$. Again
 \refT{unique} shows the claim. Finally we saw above that for
these kind of cocycles the linear form $\phi$ has to vanish on $[\ga,\ga]$.
\end{proof}

Note in particular that for $\ga$  semisimple
there will be no mixing term in \refE{sum3a} and \refE{sum3b},
see \refP{semi}.
%%%%%%%%%%%%%%%%%%%%%%%%%%%%%%%%%%%%%%%%%%%%%%%%
\subsection{Cocycles for the semisimple case}

If $\ga$ is reductive, we will sharpen the results.
First we deal with the semisimple case.
\begin{theorem}\label{T:dsimple}
(a) Let $\ga$ be a semisimple Lie algebra and $\g$ a local cocycle
of $\D_\ga$.
Then there exists a symmetric invariant bilinear form $\a$ for 
$\ga$ such that $\g$ is cohomologous to a linear combination of
the local cocycle $\g_{\a,S}$ 
given by \refE{ssimple}  and of the local
cocycle $\g_{S,R}^{(v)}$  
\refE{vecg} for $C=C_S$ of the vector field algebra $\L$.
\newline
(b) If $\ga$ is a simple Lie algebra, then $\g_{\a,S}$ is
a multiple of the standard cocycle \refE{simple} for $\ga$.
\newline
(c) $\dim\H^2_{loc}(\D_\ga,\C)=M+1$, where
$M$ is the number of simple summands of $\ga$.
\end{theorem}
\begin{proof}
Let $\g$ be a local cocycle.
By restricting it to $\L$ and $\gb$ we obtain local cocycles
 $\g_1$ and $\g_2$. By \refT{ss}, $\g_2$ is cohomologous to 
the cocycle $\g_{\a,S}$ on $\gb$ with a suitable
$\a$. The necessary coboundary $\delta\phi$
can be extended to
the whole algebra by setting the linear form $\phi$ zero on $\L$.
Hence by replacing $\g$ by a cohomologous cocycle
we  might even assume that $\g_2$ is already $\g_{\a,S}$.
But this  cocycle can be extended by zero
to the whole algebra.
Hence, $\g-(\g_1+\g_2)$ is again a cocycle
which is now of pure mixing type.
By \refP{semi} for $\ga$ semisimple they are vanishing.
This shows (a). Statement (b) follows from the uniqueness of
the invariant symmetric bilinear forms for simple Lie algebras.
Statement (c) follows from \refT{ss}.
\end{proof}
If the cocycle $\g_2$ appearing in the proof of 
\refT{dsimple} restricted to $\gb$  is $\L$-invariant
then by  \refT{ss} (b) the
cocycle $\g_2$ equals already $\g_{\a,S}$.
%%%%%%%%%%%%%%%%%%%%%%%%%%%%%%%%%%%%%%%%%%%%%%%%%%%%%%%%%%%%%%
\subsection{Cocycles for the abelian case}
%%%%%%%%%%%%%%%%%%%%%%%%%%%%%%%%%%%%%%%
Let $\ga$ be an abelian Lie algebra and $\g$ a local cocycle
for $\D_{\ga}$.
From \refE{ceAB} it follows that
\begin{equation}
\g(x(e\ldot g),y(h))+
\g(x(e\ldot g),y(h))
=0.
\end{equation}
By definition this implies that $\g_2:=\g_{|\gb}$ is
$\L$-invariant.
Also $\g_1:=\g_{|\L}$ is  a vector field cocycle.
Both cocycles can be extended by
zero on the complement to $\D_\ga$.
Hence $\g-(\g_1+\g_2)$ will be a local cocycle 
for $\D_\ga$  which is of pure 
mixing type.
\begin{proposition}
\label{P:mixd}
Every local cocycle $\g$  of pure mixing type
is cohomologous to  
\begin{equation}\label{E:amix}
\g_{\phi}(e,x(g)):=\frac {\phi(x)}{2\pi\i}
\int_{C_S}
\left(\tilde e\cdot g''+T\cdot (\tilde e\cdot
 g')\right)dz,
\end{equation}
with a suitable 
linear form
$\phi\in\ga^*$.
Vice versa, every $\phi$ defines via \refE{amix}
a local cocycle of pure mixing type.
\end{proposition}
\begin{proof}
If $\dim\ga=1$ then by a result of \cite{SchlCo}, recalled in
\refT{unique}, $\g$ is up to coboundary a scalar multiple of
$\g_{S,T}^{(m)}$, the standard mixing cocycle, see Equation
\refE{mixg}.
Let $x_1,\ldots, x_n$ be a basis of $\ga$. Set 
$\A_i:=x_i\otimes \A\cong \A$ and 
$\Da_i:=\A_i\oplus\L$.
The space $\Da_i$ is a subalgebra of $\D_\ga$ isomorphic to 
$\D$.
Restricting the cocycle to $\Da_i$ we obtain
a cocycle $\g_i$ cohomologous to 
$r_i\g_{S,T}^{(m)}$ with suitable 
$r_i\in\C$. We set $\phi(x_i):=r_i$ and
extend $\phi$ linearly to $\ga$.
If $x=\sum_i s_ix_i$ then
\begin{equation}
\g(e,x(g))=\sum_i s_i\g(e,x_i(g))\thicksim
(\sum_i s_ir_i)\g_{S,T}^{(m)}(e,g)
=\phi(x)\g_{S,T}^{(m)}(e,g).
\end{equation}
Here $\thicksim$ denotes cohomologous equivalent.
For this step the following remark has to be taken in account.
The individual cocycles $\g_i$ are cohomologous to the
standard cocycle. The individual coboundary is determined by a
linear form on $\A_i\oplus\L$.
By the pure mixing type using \refL{cobmix} the 
corresponding form can by chosen to vanish 
on the summand $\L$.
Hence after such a modification they glue together to 
a linear form on $\gb$ and define a coboundary for
$\D_\g$.
This shows the formula \refE{amix}.
Clearly, every such expression is a local cocycle of
pure mixing type.
\end{proof}
\begin{theorem}
\label{T:dabel}
Let $\ga$ be an abelian finite-dimensional Lie algebra.
The subspace of local cohomology classes
$\H^2_{loc}(\D_\ga,\C)$ of 
$\H^2(\D_\ga,\C)$ is
$\dfrac {n(n+1)}{2}+n+1$-dimensional.
Up to coboundary every local cocycle is a linear combination
of 
\begin{equation}
\g_\a(x(f),y(g))=\frac {\a(x,y)}{2\pi\i}\int_{C_S}fdg,
\end{equation}
with an arbitrary symmetric bilinear form $\a$, of
$\g_{\phi}$ with an arbitrary linear form $\phi$ (see \refE{amix}),  and of
$\g_{S,R}^{(v)}$ (see \refE{vecg}).
\end{theorem}
\begin{proof}
Let $\g$ be a local cocycle for $\D_\ga$.
As explained in the beginning of this subsection
such $\g$ can be written by restriction as
$\g=\g_{|\L}+\g_{|\gb}+\g_3$ with $\g_3$ a 
local cocycle of pure mixing type.
By \refP{mixd}, $\g_3$ can be given up to coboundary as
$\g_\phi$ with a suitable $\phi$.
The cocycle $\g_{|\L}$ is a local vector field cocycle.
Hence, $\g_{|\L}$ is a multiple of $\g_{S,R}^{(v)}$
(see \refT{unique}).
By \refT{abel}, $\g_{|\gb}=\g_\a$ with a suitable symmetric
bilinear form $\a$.
Hence every local cocycle is cohomologous to a linear
cocycle of the required type.
Vice versa, every such linear combination is a local
cocycle  and 
all the basis cocycles remain linearly independent.
\end{proof}

%%%%%%%%%%%%%%%%%%%%%%%%%%%%%%%%%%%%%%%
%%%%%%%%%%%%%%%%%%%%%%%%%%%%%%%%%%%%%%%
%%%%%%%%%%%%%%%%%%%%%%%%%%%%%%%%%%%%%%%%%%%%%%%%%%%%%%%%%%%%%%
\subsection{Cocycles for the reductive  case}

Here we have to combine the results on the semisimple and
the abelian case.
Let $\ga=\ga_0\oplus\ga_1\oplus \cdots \oplus \ga_M$ be the usual
decomposition with $\ga_0$ abelian of dimension $n$ and the 
$\ga_1,\ldots,\ga_M$ simple.
Let $\ga'=[\ga,\ga]=\ga_1\oplus \cdots \oplus \ga_M$ be the derived
subalgebra, which is now semisimple.
The algebras $\D_{\ga_0}$ and $\D_{\ga'}$ are in a natural way
subalgebras of $\D_\ga$. They are even ideals, but not
complementary ones.
\begin{proposition}\label{P:restr}
Let $\g$ be a cocycle of $\Do_{\ga}$
then $\g$ restricted to $\gb$
is $\L$-invariant if and only if $\g$ restricted to 
 $\gb'$ is  $\L$-invariant.
\end{proposition}
\begin{proof}
We only have to show that a cocycle $\g$ of 
 $\D_{\ga}$ which restricted to 
 $\gb'$
is $\L$-invariant is also  $\L$-invariant if restricted to 
 $\gb$.
By the cocycle condition \refE{ceAB} the $\L$-invariance 
on  $\gb$ is
true if and only if 
$\g(e,[x,y](gh))=0$ for $x,y\in\ga$.
$\ga$ is reductive, hence $\ga'=[\ga',\ga']$.
This implies $[x,y]=[x',y']$ with suitable $x',y'\in \ga'$,
and $\g(e,[x',y'](gh))=0$ by assumption.
\end{proof}
\begin{proposition}\label{P:linvg}
(a) Every local cocycle for $\Do_{\ga}$  restricted
to $\gb_0$ is $\L$-invariant.
\newline
(b) Every local cocycle for $\Do_{\ga}$ is cohomologous to
a local cocycle which restricted to $\ga$ is an
$\L$-invariant cocycle.
\end{proposition}
\begin{proof}
(a) Let $\g$ be a local cocycle for $\D_\ga$. It defines by restriction
to $\D_{\ga_0}$ a local cocycle.
Using \refE{ceAB} for $\D_{\ga_0}$ we obtain, 
using $[x_0,y_0]=0$, that the restriction to $\gb_0$ is
$\L$-invariant.
\newline
(b)
Let $\g$ be a local cocycle. By \refT{dsimple}, $\gamma$
restricted to $\Do_{\ga'}$
is cohomologous to a linear combination of
the $\L$-invariant cocycle \refE{ssimple} and the
separating vector field cocycle.
Let $\psi=\delta \varphi$ be the   coboundary, which appears as
difference.
The linear form $\varphi$ can be extended
to $\tilde\varphi$ by setting it zero on $\ga_0$.
We  obtain
a coboundary $\tilde\psi=\delta\tilde\varphi$
on  $\Do_{\ga}$.
The cocycle $\g'=\g+\tilde\psi$ is cohomologous to the one
we started with and restricted to $\overline{\ga'}$ it is
$\L$-invariant.
The claim (b) follows now from 
\refP{restr}.
\end{proof}
\begin{theorem}\label{T:dgred}
(a) For every local cocycle $\g$ for $\D_\ga$ there exists a symmetric
invariant bilinear form $\a$ of $\ga$, and a linear form $\phi$ of $\ga$ which
vanishes on $\ga'$, such that $\g$ is cohomologous to
\begin{equation}\label{E:gred}
\g'=r_1\g_{\a,S}+r_2\g_{\phi,S}+r_3\g_{S,R}^{(v)},
\end{equation}
with $r_1,r_2,r_3\in\C$,  a current algebra cocycle 
$\g_{\a,S}$ given by \refE{ssimple},  a mixing cocycle
$\g_{\phi,S}$ given by 
\refE{amix}, and the vector field cocycle 
$\g_{S,R}^{(v)}$ given by \refE{vecg}.
Vice versa, any such $\a$, $\phi$, $r_1,r_2,r_3\in\C$ determine a
local cocycle.
\newline
(b) 
The space of local cocycles $\H_{loc}^2(\D_\ga,\C)$ is
$\dfrac {n(n+1)}{2}+n+M+1$ dimensional.
\end{theorem}
\begin{proof}
By \refP{linvg} we may replace $\g$
by a cohomologous cocycle which restricted to  $\gb$
is $\L$-invariant. We denote it by the same symbol.
Hence by \refT{red},  $\g$ restricted to  $\gb$
is given as $\g_{\a,S}$ with $\a$ a suitable
 symmetric invariant  bilinear form.
Also by restriction to $\L$ we obtain
$\g_{|\L}=r_3\g_{S,R}^{(v)}+\psi$ with $\psi$ a coboundary of $\L$.
By the $\L$-invariance 
of $\g_{\a,S}$
all these cocycles can be extended
by zero to the whole of $\Do_{\ga}$
and we can consider the cocycle
\begin{equation}
\g'=\g-(\g_{\a,S}+
r_3\g_{S,R}^{(v)}+\psi).
\end{equation}
It will vanish on pairs of elements of the same type.
From \refE{ceAB} follows that $\g'(e,[x,y](g))=0$.
Hence by decomposing $x\in\ga$ as $x=x_0+x'$ with $x_0\in\ga_0$ and
$x'\in\ga'$ we get  
$\g'(e,x(g))=\g'(e,x_0(g))$.
It follows that $\g'$ is a local cocycle of $\D_{\ga_0}$ of 
pure mixing type which is extended by zero on $\ga'$.
In particular, using
\refT{dabel}, up to coboundary it can be given as 
$\g_{\phi,S}$ with $\phi\in\ga^*$ satisfying $\phi_{|\ga'}\equiv 0$.
This shows the first claim.
Clearly, we obtain 
a local cocycle by any such  linear combination.
Altogether (using the results  of \refT{dsimple} and \refT{dabel})
 we obtain the 
statement (b) about the dimension of the cohomology space.
\end{proof}
%%%%%%%%%%%%%%%%%%%%%%%%%%%%%%%%%%%%%%%
%%%%%%%%%%%%%%%%%%%%%%%%%%%%%%%%%%%%%%%
A further consequence is:
\begin{proposition}\label{P:gmix}
(a) 
Let $\g$ be a local cocycle of $\Do_{\ga}$.
There exists a linear form  $\phi$ on $\ga$ which vanishes
on $\ga'$ and a linear form $\psi$ on $\gb$ 
such that 
\begin{equation}\label{E:gm1}
\g(e,x(g))=\phi(x)\g^{(m)}(e,g)+\psi(x(e.g)),\quad
\forall e\in \L, g\in \A, x\in \ga,
\end{equation}
with $\g^{(m)}$ a mixing cocycle of $\Do$.
\newline
(b)
Assume in addition that the cocycle $\g$ is $\L$-invariant if restricted
to $\gb$ (or equivalently to $\gb'$) then the 
linear form $\psi$ on $\gb$ can be chosen such that it
vanishes on $\gb'$. 
\end{proposition}
\begin{proof}
(a) We  use \refT{dgred}.
If we evaluate $\g$ for pairs $(e,x(g))$ we obtain from \refE{gred}
the expression \refE{gm1}. Note that every coboundary 
will be given by a linear form $\psi$ on $\Do_\ga$, evaluated at the
Lie bracket of the two elements. But this reduces exactly to the
expression above.
\newline
(b) Now assume $\g$ to be $\L$-invariant if restricted to $\gb$.
Let $\g'$ be the corresponding cohomologous cocycle of
type \refE{gred}.
The difference will be a coboundary which evaluated at
the relevant pairs will be $\delta\psi$.
In particular it will be an $\L$-invariant cocycle.
Assume that $\psi$ does not vanish on some element of the type
$x(e\ldot g)$ with $x\in\ga'$, then $x=[y,z]$, $y,z\in\ga'$.
Hence $\psi([y,z](e\ldot g))\ne 0$ in contradiction
to the assumed $\L$-invariance, see \refE{ceAB}.
Hence, by setting $\tilde\psi(x(h))=0$ for all   
$x\in\ga'$, $h\in\A$ we obtain
$\psi(x(h))=\tilde\psi(x(h))$ and get the required form
in \refE{gm1}.
\end{proof}
%%%%%%%%%%%%%%%%%%%%%%%%%%%%%%%
\section{Examples: $\sln(n)$ and $\gl(n)$}
\label{S:special}
%\input example.tex
%%%%%  Examples     17.10.02
%%%%%%%%%%%%%%%%%%%%%%%%%%%%%%%
As examples we deal with the important special cases 
 $\sln(n)$ (which is a simple algebra) and the case of 
 $\gl(n)$ (which is reductive but not semi-simple).
We will study the affine algebras, the differential operator
algebras, and their central extensions.
The corresponding algebras appear in the study of
fermionic representations, see \cite{Shferm}.

\subsection{ $\sln(n)$}
%%%%%%%%%%%%%%%%%%%%%%%%%%%%%%%%%%%%%
First we consider $\sln(n)$, the Lie algebra of trace-less complex $n\times n$
matrices.
Up to multiplication with a scalar 
the Cartan-Killing form $\b(x,y)=\tr(xy)$ is the 
unique
symmetric invariant bilinear form.
It is non-degenerate.

Due to the fact that $\sln(n)$ is simple, from 
\refT{simple} and \refT{dsimple}
follows
\begin{proposition}\label{P:slcoc}
(a) Every local cocycle for the current algebra $\overline{\sln}(n)$
is cohomologous to
\begin{equation}\label{E:slcoc}
\g(x(g),y(h))=r
\frac {\tr(xy)}{2\pi\i}\int_{C_S} gdh,\quad r\in\C.
\end{equation}
(b) Every $\L$-invariant local cocycle equals the cocycle \refE{slcoc}
with a suitable $r$.
\newline
(c) Every local cocycle for the differential operator algebra
$\Do_{\sln(n)}$ is cohomologous to a linear combination of
\refE{slcoc} and the standard local cocycle $\g_{S,R}^{(v)}$ 
for the vector field
algebra.
In particular, there exist no cocycles of pure mixing type.
\end{proposition}
%%%%%%%%%%%%%%%%%%%%%
\subsection{ $\gl(n)$}
Next we deal with $\gl(n)$, the Lie algebra of all complex
$n\times n$-matrices.
Recall that $\gl(n)$ can be written  as Lie algebra direct sum
$\gl(n)=\mathfrak{s}(n)\oplus \sln(n)\cong\C \oplus\sln(n)$.
Here $\mathfrak{s}(n)$ denotes the $n\times n$ scalar matrices.
This decomposition is the decomposition as reductive Lie 
algebra into its
abelian and semi-simple summand.

After tensoring with $\A$ we obtain 
$\overline{\gl}(n)=\overline{\mathfrak{s}}(n)\oplus \overline{\sln}(n)\oplus 
\cong \A\oplus\overline{\sln}(n)$.
From \refT{red} it follows that we need to determine all
symmetric invariant bilinear forms for $\overline{\gl}(n)$.
From the above decomposition it follows that the space of
such forms  is
two-dimensional.
An adapted basis is
\begin{equation}
\a_1(x,y)=\tr(xy),\quad\text{and}\quad
\a_2(x,y)=\tr(x)\tr(y).
\end{equation}
The form $\a_1$ is the ``natural'' extension for the Cartan-Killing
form for $\sln(n)$ to $\gl(n)$ and is also $\gl(n)$ invariant.
The form $\a_2$ vanishes on $\sln(n)$ and is non-zero on $\mathfrak{s}(n)$.
Note that ${\a_1}_{|\mathfrak{s}(n)}\not\equiv 0$ .
The by zero on $\mathfrak{s}(n)$ to $\gl(n)$  extended form 
$\tilde\a_1$ of the Cartan-Killing form 
for $\sln(n)$ would be
\begin{equation}
\tilde\a_1(x,y)=\tr(xy)-\frac 1n\tr(x)\tr(y).
\end{equation}
%%%%%%%%%%%%%%%%%%%%%%%%%%%%%%%%%%%%%%%%%
From \refT{red} we conclude
%%%%%%%%%%%%%%%%%%%
\begin{proposition}\label{P:glcoc}
(a) A cocycle for $\overline{\gl}(n)$ is local and 
restricted to $\overline{\mathfrak{s}}$ is $\L$-invariant
if and only if it is cohomologous to a linear combination of the following
two cocycles
\begin{equation}\label{E:glns}
\g_1(x(g),y(h))=\frac{\tr(xy)}{2\pi\i} \int_{C_S} gdh,\qquad
\g_2(x(g),y(h))=\frac{\tr(x)\tr(y)}{2\pi\i} \int_{C_S} gdh.
\end{equation}
(b) If the cocycle $\g$ of (a) restricted to 
$\overline{\gl}(n)$ is $\L$-invariant then 
$\g$ is equal to the linear combination of the cocycles
\refE{glns}
\end{proposition}
\begin{proposition}\label{P:dgln}
(a) Every local cocycle $\g$ for  $\Do_{\gl(n)}$ is cohomologous to a
linear combination of the cocycles $\g_1$ and $\g_2$ of \refE{glns},
of the mixing  cocycle
\begin{equation}
\g_{3,T}(e,x(g))=\frac {\tr(x)}{2\pi\i}\int_{C_S}
\big(\tilde e g''+T\tilde e g'\big)dz,
\end{equation}
and of the standard local cocycle $\g_{S,R}^{(v)}$ 
for the vector field algebra, i.e.
\begin{equation}\label{E:decompc}
\gamma=r_1\gamma_1+r_2\gamma_2
+r_3\g_{3,T}+r_4\g_{S,R}^{(v)}+\text{coboundary},
\end{equation}
with suitable $r_1,r_2,r_3,r_4\in\C$.
\newline
(b) If the cocycle $\g$ is local and 
restricted to $\overline{\gl}(n)$ is $\L$-invariant, 
and
$r_3,r_4\ne 0$ then there exist an affine connection $T'$
and a projective connection $R'$ holomorphic
outside $A$ such that 
$\gamma=r_1\gamma_1+r_2\gamma_2
+r_3\g_{3,T'}+r_4\g_{S,R'}^{(v)}$.
\newline
(c) $\dim\H^2_{loc}(\Do_{\gl(n)},\C)=4$.
\end{proposition}
\begin{proof}
For $\gl(n)=\mathfrak{s}(n)\oplus\sln(n)$ up to multiplication
with a scalar, the unique linear form vanishing on $\sln(n)$ is
$x\mapsto\tr(x)$.
Hence from \refP{glcoc} and \refT{dgred} the claim (a) follows.
To prove (b) we have to note that by the
$\L$-invariance
using \refP{glcoc}(b) it follows that
the restriction to $\overline{\gl}(n)$ is without
changing the element in the cohomology class  
already of the required type. 
Also by the non-vanishing  of the coefficients $r_3$ and
$r_4$ the coboundaries can be included into the definition of 
$\g_{S,R'}^{(v)}$ and
$\g_{3,T'}$ by choosing suitable projective and affine projections.
Part (c) follows from \refT{dgred}, Part (b).
\end{proof}
%%%%%%%%%%%%%%%%%%%%%%%%%%%%%%%%%%%%%%%%%%
For the cocycle obtained by the regularization 
 of the fermionic representations associated to  $\gl(n)$
its mixing part  was studied by Sheinman in \cite{Shsc}.
Independent of the results in \cite{SchlCo},
he obtains results on its structure.
%%%%%%%%%%%%%%%%%%%%%%%%%%
\appendix
\section{Proof of \refL{split}}
\label{S:split}
%\input appendix.tex
%   Appendix      
%%%%   16.10.02
%%%%%%%%%%%
In this appendix I give the details of the proof of 
\refL{split} by recalling the steps of the proof of the
uniqueness for the local and $\L$-invariant cocycles for the 
function algebra $\A$ (see \cite[Sections 5.1 and 5.2]{SchlCo})
now suitable adjusted to our more general situation.
Let $\g$ be an $\L$-invariant and local cocycle for the
current algebra $\gb=\ga\otimes\A$ of an abelian finite-dimensional
Lie algebra $\ga$.
We fix $x,y\in\ga$ and set $\g_{x,y}(f,g):=\g(x\otimes f,y\otimes g)$.
This will define a bilinear form on $\A$. A priori we do not
know that $\gxy$ will be antisymmetric.
By the  $\L$-invariance we obtain
\begin{equation}
\label{E:alinv}
\gxy(e\ldot f,g)+\gxy(f,e\ldot g)=0,\quad
\forall e\in\L,\ \forall f,g\in\A.
\end{equation}
Our aim is to show that we can choose $\a_{xy}\in\C$ such that
$\gxy(A_{m,s},A_{l,t})=\a_{xy}\g_S^{(f)}(A_{m,s},A_{l,t})$
for all pairs of the basis elements 
$\{A_{m,s}\mid m\in\Z,\ s=1,\ldots,K\}$ of $\A$.
  
By the almost-graded action of $\L$ on $\A$ we have
\begin{equation}\label{E:alf}
e_{n,r}.A_{m,s}=\delta_r^s \cdot m\cdot A_{n+m,r}
+\sum_{h=n+m+1}^{n+m+L_2}\sum_{t=1}^K 
 b_{(n,r),(m,s)}^{(h,t)} A_{h,t},
\end{equation}
with $b_{(n,r),(m,s)}^{(h,t)}\in\C$, and $L_2$ the upper bound
for the almost-graded structure.
Equation  \refE{alinv} for  basis elements writes as
\begin{equation}
\gxy(e_{n,r}\ldot A_{m,s},A_{p,t})+
\gxy( A_{m,s}, e_{n,r}\ldot A_{p,t})=0.
\end{equation}
Using \refE{alf} we get
\begin{equation}\label{E:induct}
\delta_r^s\cdot m\cdot \gxy(A_{m+n,r},A_{p,t})+
\delta_r^t\cdot p\cdot \gxy(A_{m,s},A_{n+p,t})=\hL
\end{equation}
Here and in the following I will  use the phrase
``can be expressed by elements of higher level'',
``determined by higher level'', or simply 
``$=\hL$'' to denote that it is a universal linear combination
of  values for pairs of homogeneous
elements of levels higher than the level under consideration.
By the level I understand the sum of the degrees of the two
arguments.
The coefficients appearing in this linear combination only
depend on the  geometric situation, i.e. on the 
structure constants of the algebra  and not on the
the bilinear form  under consideration.
In particular, if two bilinear forms 
are given by higher level and they coincide in  higher levels, they will
coincide also for the elements under consideration.

For $r= t\ne s$ we obtain
\begin{equation}
p\cdot\gxy(A_{m,s},A_{n+p,t})=\hL
\end{equation}
Hence, as long as $s\ne t$ the  value
$\gxy(A_{m,s},A_{l,t})$ is determined by higher level.
Next we 
consider $s=t=r$.
From \refE{induct} it follows 
\begin{equation}\label{E:induct2}
m\cdot\gxy(A_{m+n,r},A_{p,r})+p\cdot\gxy(A_{m,r},A_{n+p,r})=\hL
\end{equation}
We set $n=0$ and obtain
\begin{equation}
(m+p)\cdot\gxy(A_{m,r},A_{p,r})=
\hL
\end{equation}
Hence, for $p\ne -m$ the  value
$\gxy(A_{m,r},A_{p,r})$ is determined by higher level.

Next in \refE{induct2} we set  $p=-(n+1), m=1$ and obtain
\begin{equation}
\gxy(A_{n+1,r},A_{-(n+1),r})=
(n+1)\gxy(A_{1,r},A_{-1,r})+\hL
\end{equation}
Now we use that $\g$ is local.
In particular, it is bounded from above.
Starting from a certain level the value of the bilinear form 
will be zero if evaluated for elements of level higher than
this level.
But by the above inductive procedure it has to vanish for all
levels greater zero.

Let $C_i$ be a (possibly deformed) circle around the incoming point $P_i\in I$.
We set
$\g_i^{(f)}(f,g):=\frac {1}{2\pi\i}\int_{C_i} fdg$.
This is a cocycle for $\A$ which is bounded from above by zero.
Obviously, for the separating cocycle
$\g_S^{(f)}=\sum_{i=1}^K\g_i^{(f)}$.
Note that $\g_i^{(f)}$ is not necessarily bounded from below and
hence not necessarily local.
By calculations of residues we obtain 
\begin{equation}
\g_i^{(f)}(A_{n,r},A_{-n,r})=-n\cdot\delta_r^i.
\end{equation}
Set 
$\tilde\g=\sum_{r=1}^K\a_{xy}^r\g_r^{(f)}$
with $\a_{xy}^r:=-\gxy(A_{1,r},A_{-1,r})\in\C$.
At level zero (and higher) $\tilde\g=\g_{xy}$, and hence by induction
at every level.
In particular  $\gxy$ is a local cocycle and we can apply the results
of \cite{SchlCo}.
By Theorem 4.3 in  \cite{SchlCo},  $\gxy$ will be  local
if and only if in above  
linear combination all coefficients are the same, i.e.  
$\a_{xy}:=\a_{xy}^1= \a_{xy}^2=\cdots= \a_{xy}^K$.
Hence, $\gxy=\a_{xy}\g_S^{(f)}$.
This proves \refL{split}.
%%%%%%%%%%%%%%%%%%%%%%%%%%%%%%%%%%%%%%%%%%%%%%%%%%%%%%%%%%%%%%
%%%%%%%%%%%%%%%%%%%%%%%%%%%%%%%%%%%%%%%%%%%%%%%%%%%%%%%%%%%%%%
%%%%%%%%%%  bibliography  %%%%%%%%%%
%\bibliographystyle{amsplain}
%\bibliography{database}
%\input affmain.bbl
%\providecommand{\bysame}{\leavevmode\hbox to3em{\hrulefill}\thinspace}
%\providecommand{\MR}{\relax\ifhmode\unskip\space\fi MR }
%\providecommand{\MRhref}[2]{%
%  \href{http://www.ams.org/mathscinet-getitem?mr=#1}{#2}
%}
%\providecommand{\href}[2]{#2}

%%%%%%%%%%%%%%%%%%%%%%%%%%%%%%%%%%%%%%%%%%%%%%%%%%%%
\end{document}